\documentclass{art_ar}
\usepackage{graphicx}
\usepackage{amssymb}

\usepackage{pstricks, pst-plot, pst-node, pst-tree}
\usepackage{array}
\usepackage{amsmath}
\usepackage{enumerate}
\usepackage{longtable}
\usepackage{rotating}

\theoremstyle{plain}
\newtheorem{Def}{Definition}[section]
\newtheorem{Sat}[Def]{Proposition}
\newtheorem{The}[Def]{Theorem}
\newtheorem{Kor}[Def]{Corollary}

\newcommand{\E}{\operatorname{E}}
\newcommand{\Prob}{\operatorname{P}}

\begin{document}

\begin{frontmatter}



\title{Stochastic Taylor Expansions for Functionals of Diffusion Processes}


\author{Andreas R{\"o}{\ss}ler}

\address{Technische~Universit\"at~Darmstadt, Fachbereich~Mathematik, Schlossgartenstr.~7,
D-64289~Darmstadt, Germany}

\ead{roessler@mathematik.tu-darmstadt.de}

%\ead[url]{http://www.mathematik.tu-darmstadt.de/\symbol{126}roessler/}

\begin{abstract}
% Text of abstract
In the present paper, a stochastic Taylor expansion of some
functional applied to the solution process of an It{\^o} or
Stratonovich stochastic differential equation with a
multi-dimensional driving Wiener process is given. Therefore, the
multi-colored rooted tree analysis is applied in order to obtain a
transparent representation of the expansion which is similar to
the B-series expansion for solutions of ordinary differential
equations in the deterministic setting. Further, some estimates
for the mean--square and the mean truncation errors are given.
\end{abstract}

\begin{keyword}
stochastic Taylor expansion \sep stochastic differential equation
\sep multi--colored rooted tree analysis \sep strong approximation
\sep mean-square approximation
%\PACS
%\MCS
\\ {\emph{MSC:}} 60H35 \sep 41A58 \sep 65C30 \sep 60H10 \sep 60J60
\end{keyword}
\end{frontmatter}
%
% main text
\section{Introduction}
\label{intro}
Let $(\Omega, \mathcal{F}, \Prob)$ be a probability space with a
filtration $(\mathcal{F}_t)_{t \geq 0}$ fulfilling the usual
conditions.
% and let $\mathcal{I}=[t_0,T]$ for some $0 \leq t_0 < T < \infty$.
Since each non-autonomous stochastic differential equation (SDE)
can be written as an autonomous SDE system with one additional
equation representing time, we consider without loss of generality
autonomous SDE systems only. Thus, for some $0 \leq t_0 < T <
\infty$ let $(X_t)_{t \in [t_0,T]}$ be the solution of the
$d$--dimensional autonomous SDE system
\begin{equation} \label{St-lg-sde-ito-strato-1}
    {\mathrm{d}} X_t = a(X_t) \, {\mathrm{d}}t + b(X_t) \,
    *{\mathrm{d}}W_t
\end{equation}
with an $m$--dimensional driving Wiener process $(W_t)_{t \geq 0}$
% =((W_t^1, \ldots, W_t^m))_{t \geq 0}$
w.r.t.\ $(\mathcal{F}_t)_{t \geq 0}$. Then,
SDE~(\ref{St-lg-sde-ito-strato-1}) can be written in integral form
\begin{equation} \label{Intro-Ito-St-SDE1-integralform-Wm}
    X_t = X_{t_0} + \int_{t_0}^t a(X_s) \, {\mathrm{d}}s + \sum_{j=1}^m
    \int_{t_0}^t b^j(X_s) \, *{\mathrm{d}}W_s^j
\end{equation}
for $d,m \geq 1$ and $t \in [t_0,T]$, where we write
$*{\mathrm{d}}W_s^j = {\mathrm{d}}W_s^j$ in the case of an It{\^o}
stochastic integral and $*{\mathrm{d}}W_s^j = \circ
{\mathrm{d}}W_s^j$ for a Stratonovich stochastic integral. Here,
we suppose that $a: \mathbb{R}^d \to \mathbb{R}^d$ and $b :
\mathbb{R}^d \to \mathbb{R}^{d \times m}$ are measureable
functions which are sufficiently smooth and we denote by $b^j$ for
$j=1, \ldots,m$ the $j$th column of the $d \times m$-matrix
function $b=(b^{i,j})$.
Let $X_{t_0} \in \mathbb{R}^d$ be the
$\mathcal{F}_{t_0}$-measurable initial value with $X_{t_0} \in
L^2(\Omega)$.
%for some $l \in \mathbb{N}$ holds $\E(\|X_{t_0}\|^{2l})<\infty$
In the following, we suppose that the conditions of the Existence
and Uniqueness Theorem~\cite{KS99} are fulfilled for
SDE~(\ref{Intro-Ito-St-SDE1-integralform-Wm}) and we denote by $\|
\cdot \|$ the Euclidean norm.

The aim of the present paper is to give an expansion of $f(X_t)$
for some functional $f : \mathbb{R}^d \to \mathbb{R}$. Therefore,
we define for $j=1, \ldots, m$ the operators
\begin{equation}
    \hat{L}^{0} = \sum_{k=1}^d a^k \frac{\partial}{\partial x^k},
    \quad \quad
    \hat{L}^{j} = \sum_{k,l=1}^d b^{k,j} \, b^{l,j}
    \frac{\partial^2}{\partial x^k \partial x^l}, \quad \quad
    L^j = \sum_{k=1}^d b^{k,j} \frac{\partial}{\partial x^k},
\end{equation}
and $L^0 = \hat{L}^{0} + \tfrac{1}{2} \sum_{j=1}^m \hat{L}^{j}$.
Considering now the It{\^o}
SDE~(\ref{Intro-Ito-St-SDE1-integralform-Wm}), we obtain for
sufficiently smooth $f$ by recursive application of It{\^o}'s
formula
\begin{equation} \label{ito-taylor-expansion-bsp1}
    \begin{split}
    f(X_t) &= f(X_{t_0}) + \int_{t_0}^t L^0 f(X_s) \, {\mathrm{d}}s
    + \sum_{j=1}^m \int_{t_0}^t L^j f(X_s) \, {\mathrm{d}}W_s^j \\
    &= f(X_{t_0}) + \int_{t_0}^t L^0 f(X_s) \, {\mathrm{d}}s
    + \sum_{j_1=1}^m \int_{t_0}^t  \left( L^{j_1} f(X_{t_0})
    + \int_{t_0}^s L^0 L^{j_1} f(X_u) \, {\mathrm{d}}u \right. \\
    &+ \sum_{j_2=1}^m \left. \int_{t_0}^s L^{j_2} L^{j_1} f(X_u)
    \, {\mathrm{d}}W_u^{j_2} \right) \, {\mathrm{d}}W_s^{j_1} \, .
    \end{split}
\end{equation}
Repetition of this procedure by applying It{\^o}'s formula to $L^0
f(X_s)$ yields a further expansion, and so on. As a result of
this, we obtain the It{\^o}-Taylor-expansion due to Platen and
Wagner~\cite{KP99,PlWa82} with a remainder of integral type. In
the following, we always assume that all necessary derivatives and
multiple stochastic integrals exist. If
SDE~(\ref{Intro-Ito-St-SDE1-integralform-Wm}) is considered in the
Stratonovich sense, then the same expansion as in
(\ref{ito-taylor-expansion-bsp1}) applies however with $L^0$
replaced by $\hat{L}^0$.

In the present paper, we develop an expansion of $f(X_t)$ based on
multi--colored rooted trees. This turns out to be an extension of
the rooted tree approach for ordinary differential equations due
to Butcher~\cite{Butcher03} in the deterministic setting. Burrage
and Burrage~\cite{BuBu00a} developed the expansion of $f(X_t)$ by
multi--colored rooted trees for the special case when $X_t$ is the
solution of a Stratonovich SDE and when $f(X_t)=X_t$. In contrast
to this, we give an expansion not only for the Stratonovich
solution $X_t$ but also for solutions of It{\^o} SDEs and additionally
for arbitrary functionals $f(X_t)$ of the solution process.
Therefore, we follow the approach proposed in \cite{NNRT07,Roe04b}
and make use of an additional node corresponding to the functional
$f$ for the root of the considered trees. Further, we have to take
into account the more complex structure of the operator $L^0$ for
It{\^o} SDEs compared to $\hat{L}^0$ for Stratonovich SDEs. As the
main advantage of the rooted tree expansion of $f(X_t)$, we obtain
a clear and simple expansion with equal elementary differentials
pooled together. Compared to the approach based on hierarchical
sets by Kloeden and Platen~\cite{KP99}, each elementary
differential can be determined directly by the corresponding
rooted tree. Further, expansions based on rooted trees allow a
systematic development of higher order derivative free
approximations similar to the deterministic
setting~\cite{BuBu00a,Butcher03,Roe06a}.
%
%
%
% ==============================================================
\section{Colored Rooted Tree Analysis}
\label{Sec:Order-Conditions}
% ==============================================================
%
%
%
Following the approach in \cite{Roe06a,Roe04b}, we give a
definition of colored trees which will be suitable for SDEs
w.r.t.\ a multi--dimensional Wiener process.
%
%Since each SDE system can be represented by an autonomous SDE
%system
%%
%\begin{equation} \label{Ito-St-SDE1-autonomous-Wm}
%    X_t = X_{t_0} + \int_{t_0}^t a(X_s) \, {\mathrm{d}}s + \sum_{j=1}^m
%    \int_{t_0}^t b^j(X_s) \, *{\mathrm{d}}W_s^j
%\end{equation}
%%
%with one additional equation representing time, we restrict our
%considerations to an autonomous SDE system in this section.
%

\begin{Def} \label{Def:rooted-S-trees:Wm}
%    Let $l$ be a positive integer.
%    \begin{enumerate}
%        \item
        A monotonically labelled {\emph{S-tree (stochastic
        tree)}} $\textbf{t}$ with $l=l(\textbf{t}) \in \mathbb{N}$
        nodes is a pair of maps
        $\textbf{t}=(\textbf{t}',\textbf{t}'')$ with
        \begin{equation*}
            \begin{split}
                \textbf{t}' & : \{2, \ldots, l\} \to \{1, \ldots,
                l-1\} \\
                \textbf{t}'' & : \{1, \ldots, l\} \to \mathcal{A} \\
            \end{split}
        \end{equation*}
        so that $\textbf{t}'(i) < i$ for $i=2, \ldots, l$. Unless otherwise noted, we choose
        the set $\mathcal{A} = \{ \gamma, \tau_{j} : j \in \{0,1, \ldots, m\}
        \}$.
%        where $j_k$ is a variable index with $j_k \in \{1, \ldots, m\}$.
%        \item
        Let $LTS$ denote the set of all monotonically labelled
        S-trees w.r.t.\ $\mathcal{A}$.
%        Here two trees $\textbf{t}=(\textbf{t}',\textbf{t}'')$ and $\textbf{u}=(\textbf{u}',\textbf{u}'')$
%        just differing by their colors $\textbf{t}''$ and $\textbf{u}''$ are
%        considered to be identical if there exists a bijective
%        map $\pi : \mathcal{A} \to \mathcal{A}$
%        with $\pi(\gamma)=\gamma$ and $\pi(\tau)=\tau$ so that
%        $\textbf{t}''(i) = \pi(\textbf{u}''(i))$ holds for $i=1, \ldots,
%        l$.
%    \end{enumerate}
\end{Def}

Then $\textbf{t}'$ defines a father son relation between the nodes,
i.e.\ $\textbf{t}'(i)$ is the father of the son $i$. Furthermore the
color $\textbf{t}''(i)$, which consists of one element of the set
$\mathcal{A}$, is added to the node $i$ for $i=1, \ldots,
l(\textbf{t})$. Here, $\tau_0 = $~\pstree[treemode=U,
dotstyle=otimes, dotsize=3.2mm, levelsep=0.1cm, radius=1.6mm,
treefit=loose]
    {\Tn}{
    \pstree[treemode=U, dotstyle=otimes, dotsize=3.2mm, levelsep=0cm, radius=1.6mm, treefit=loose]
    {\TC*~[tnpos=r]{}} {}
    }
is a deterministic node, $\tau_{j} = $~\pstree[treemode=U,
dotstyle=otimes, dotsize=3.2mm, levelsep=0.1cm, radius=1.6mm,
treefit=loose]
    {\Tn}{
    \pstree[treemode=U, dotstyle=otimes, dotsize=3.2mm, levelsep=0cm, radius=1.6mm, treefit=loose]
    {\TC~[tnpos=r]{$\!\!{_{j}}$}} {}
    }
is a stochastic node with $j \in \{1, \ldots, m\}$ and $\gamma =
\,\,$~\pstree[treemode=U, dotstyle=otimes, dotsize=3.2mm,
levelsep=0.1cm, radius=1.6mm, treefit=loose]
    {\Tn}{
    \pstree[treemode=U, dotstyle=otimes, dotsize=3.2mm, levelsep=0cm, radius=1.6mm, treefit=loose]
    {\Tdot~[tnpos=r]{ }} {}
    }
can be the root of a tree. The variable index $j$ is associated
with the $j$th component of the corresponding $m$-dimensional
Wiener process of the considered SDE.
\begin{figure}[htbp]
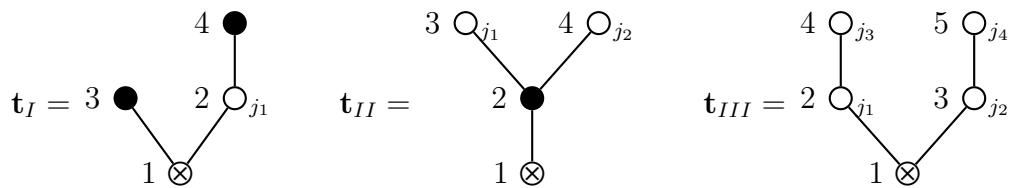

\begin{center}
\begin{tabular}{ccccc}
    $\textbf{t}_I = \begin{array}{c}
    \text{\pstree[treemode=U, dotstyle=otimes, dotsize=3.2mm, levelsep=0.1cm, radius=1.6mm, treefit=loose]
    {\Tn}{
    \pstree[treemode=U, dotstyle=otimes, dotsize=3.2mm, levelsep=1cm, radius=1.6mm, treefit=loose]
    {\Tdot~[tnpos=l]{1 }}
    {\TC*~[tnpos=l]{3}
    \pstree{\TC~[tnpos=l]{2}~[tnpos=r]{$\!\!_{j_1}$}}{\TC*~[tnpos=l]{4}}}
    }}
    \end{array}$
    & \quad &
    $\textbf{t}_{II} = \begin{array}{c}
    \text{\pstree[treemode=U, dotstyle=otimes, dotsize=3.2mm, levelsep=0.1cm, radius=1.6mm, treefit=loose]
    {\Tn}{
    \pstree[treemode=U, dotstyle=otimes, dotsize=3.2mm, levelsep=1cm, radius=1.6mm, treefit=loose]
    {\Tdot~[tnpos=l]{1 }} {\pstree{\TC*~[tnpos=l]{2}}{\TC~[tnpos=l]{3}~[tnpos=r]{$\!\!_{j_1}$}
    \TC~[tnpos=l]{4}~[tnpos=r]{$\!\!_{j_2}$}}}
    }}
    \end{array}$
    & \quad &
    $\textbf{t}_{III} = \begin{array}{c}
    \text{\pstree[treemode=U, dotstyle=otimes, dotsize=3.2mm, levelsep=0.1cm, radius=1.6mm, treefit=loose]
    {\Tn}{
    \pstree[treemode=U, dotstyle=otimes, dotsize=3.2mm, levelsep=1cm, radius=1.6mm, treefit=loose]
    {\Tdot~[tnpos=l]{1 }}
    {\pstree{\TC~[tnpos=l]{2}~[tnpos=r]{$\!\!_{j_1}$}}{\TC~[tnpos=l]{4}~[tnpos=r]{$\!\!_{j_3}$}}
    \pstree{\TC~[tnpos=l]{3}~[tnpos=r]{$\!\!_{j_2}$}}{\TC~[tnpos=l]{5}~[tnpos=r]{$\!\!_{j_4}$}}}
    }}
    \end{array}$
\end{tabular}
\caption{Three elements of $LTS$ with $j_1, j_2, j_3, j_4 \in \{1,
\ldots, m\}$.} \label{St-S-tree-examples-tI+tII:Wm}
\end{center}
\end{figure}
As an example Figure~\ref{St-S-tree-examples-tI+tII:Wm} presents
three elements of $LTS$.

In the following, we denote by $d(\textbf{t}) = |\{ i :
\textbf{t}''(i) = \tau_0 \}|$ the number of deterministic nodes
and by $s(\textbf{t}) = |\{ i : \textbf{t}''(i) = \tau_{j}, 1 \leq
j \leq m \}|$ the number of stochastic nodes. The order
$\rho(\textbf{t})$ of the tree $\textbf{t}$ is defined as
$\rho(\textbf{t}) = d(\textbf{t}) + \tfrac{1}{2} s(\textbf{t})$
with $\rho(\gamma) = 0$. The order of the trees presented in
Figure~\ref{St-S-tree-examples-tI+tII:Wm} can be calculated as
$\rho(\textbf{t}_I)=2.5$ and
$\rho(\textbf{t}_{II})=\rho(\textbf{t}_{III})=2$.

%
%\begin{figure}[htbp]
%\begin{center}
%\begin{tabular}{ccccc}
%    \pstree[treemode=U, dotstyle=otimes, dotsize=3.2mm, levelsep=0.1cm, radius=1.6mm, treefit=loose]
%    {\Tn}{
%    \pstree[treemode=U, dotstyle=otimes, dotsize=3.2mm, levelsep=1cm, radius=1.6mm, treefit=loose, nodesepB=1mm]
%    {\Tdot} {\Tr{$\textbf{t}_1$} \Tr{$\textbf{t}_2$} \Tr[edge=none]{$\cdots$} \Tr{$\textbf{t}_k$}}
%    }
%    & &
%    \pstree[treemode=U, dotstyle=otimes, dotsize=3.2mm, levelsep=0.1cm, radius=1.6mm, treefit=loose]
%    {\Tn}{
%    \pstree[treemode=U, dotstyle=otimes, dotsize=3.2mm, levelsep=1cm, radius=1.6mm, treefit=loose, nodesepB=1mm]
%    {\TC*} {\Tr{$\textbf{t}_1$} \Tr{$\textbf{t}_2$} \Tr[edge=none]{$\cdots$} \Tr{$\textbf{t}_k$}}
%    }
%    & &
%    \pstree[treemode=U, dotstyle=otimes, dotsize=3.2mm, levelsep=0.1cm, radius=1.6mm, treefit=loose]
%    {\Tn}{
%    \pstree[treemode=U, dotstyle=otimes, dotsize=3.2mm, levelsep=1cm, radius=1.6mm, treefit=loose, nodesepB=1mm]
%    {\TC~[tnpos=r]{$\!\!_{j}$}} {\Tr{$\textbf{t}_1$} \Tr{$\textbf{t}_2$} \Tr[edge=none]{$\cdots$} \Tr{$\textbf{t}_k$}}
%    }
%    \\
%    $(\textbf{t}_1, \ldots, \textbf{t}_k)$ & &
%    $[\textbf{t}_1, \ldots, \textbf{t}_k]$ & & $\,\, \{\textbf{t}_1, \ldots, \textbf{t}_k\}_j$
%\end{tabular}
%\caption{Writing a colored S-tree with brackets.}
%\label{St-tree-bracket-together}
%\end{center}
%\end{figure}
%
%
Every labelled tree can be written by a combination of brackets:
If $\textbf{t}_1, \ldots, \textbf{t}_k$ are colored trees then we
denote by $[\textbf{t}_1, \ldots, \textbf{t}_k]_{\gamma}$ and
$[\textbf{t}_1, \ldots, \textbf{t}_k]_j$ the tree in which
$\textbf{t}_1, \ldots, \textbf{t}_k$ are each joined by a single
branch to
    \mbox{$\gamma = \,\,$ \pstree[treemode=U, dotstyle=otimes, dotsize=3.2mm,
    levelsep=0.1cm, radius=1.6mm, treefit=loose]
    {\Tn}{
    \pstree[treemode=U, dotstyle=otimes, dotsize=3.2mm, levelsep=0cm, radius=1.6mm, treefit=loose]
    {\Tdot} {}
    }~$\,$}
%    \pstree[treemode=U, dotstyle=otimes, dotsize=3.2mm, levelsep=0.1cm, radius=1.6mm, treefit=loose]
%    {\Tn}{
%    \pstree[treemode=U, dotstyle=otimes, dotsize=3.2mm, levelsep=0cm, radius=1.6mm, treefit=loose]
%    {\TC*} {}
%    }~$\,\,$
    and
    $\tau_j = \,\,$\pstree[treemode=U, dotstyle=otimes, dotsize=3.2mm, levelsep=0.1cm, radius=1.6mm, treefit=loose]
    {\Tn}{
    \pstree[treemode=U, dotstyle=otimes, dotsize=3.2mm, levelsep=0cm, radius=1.6mm, treefit=loose]
    {\TC~[tnpos=r]{$\!\!_j$}} {}
    } for $j = 0,1, \ldots, m$,
    respectively. % (see Figure~\ref{St-tree-bracket-together}).
Therefore proceeding recursively, for the three examples in
Figure~\ref{St-S-tree-examples-tI+tII:Wm} we obtain
    $\textbf{t}_I =
%    [\text{\pstree[treemode=U, dotstyle=otimes, dotsize=3.2mm, levelsep=0.1cm, radius=1.6mm, treefit=loose]
%    {\Tn}{
%    \pstree[treemode=U, dotstyle=otimes, dotsize=3.2mm, levelsep=0cm, radius=1.6mm, treefit=loose]
%    {\TC*~[tnpos=r]{$\!\!^3$}} {}
%    }},
%    [\text{\pstree[treemode=U, dotstyle=otimes, dotsize=3.2mm, levelsep=0.1cm, radius=1.6mm, treefit=loose]
%    {\Tn}{
%    \pstree[treemode=U, dotstyle=otimes, dotsize=3.2mm, levelsep=0cm, radius=1.6mm, treefit=loose]
%    {\TC*~[tnpos=r]{$\!\!^4$}} {}
%    }} ]_{j_1}^2 ]_{\gamma}^1 =
    [\tau_{0}^3, [\tau_0^4 ]_{j_1}^2]_{\gamma}^1$,
    $\textbf{t}_{II} =
%    [[\text{\pstree[treemode=U, dotstyle=otimes, dotsize=3.2mm, levelsep=0.1cm, radius=1.6mm, treefit=loose]
%    {\Tn}{
%    \pstree[treemode=U, dotstyle=otimes, dotsize=3.2mm, levelsep=0cm, radius=1.6mm, treefit=loose]
%    {\TC*~[tnpos=r]{$\!\!$}} {}
%    }},
%    \text{\pstree[treemode=U, dotstyle=otimes, dotsize=3.2mm, levelsep=0.1cm, radius=1.6mm, treefit=loose]
%    {\Tn}{
%    \pstree[treemode=U, dotstyle=otimes, dotsize=3.2mm, levelsep=0cm, radius=1.6mm, treefit=loose]
%    {\TC~[tnpos=r]{$\!\!_{j_2}$}} {}
%    }}
%    ]_{j_1}]_{\gamma} =
    [[ \tau_{j_1}^3, \tau_{j_2}^4 ]_{0}^2]_{\gamma}^1$ and
    $\textbf{t}_{III} =
    [[\tau_{j_3}^4]_{j_1}^2, [\tau_{j_4}^5]_{j_2}^3]_{\gamma}^1$
for $j_1, j_2, j_3, j_4 \in \{1, \ldots, m\}$.
%
% ===========================================================
%

Now, two labelled trees $\textbf{t}, \textbf{u} \in LTS$ with
$l=l(\textbf{t})=l(\textbf{u})$ nodes are called equivalent, i.e.\
$\textbf{t} \sim \textbf{u}$, if there exists a bijective map $\pi
: \{1, \ldots, l\} \to \{1, \ldots, l\}$ with $\textbf{t}'(i) =
\pi^{-1}(\textbf{u}'(\pi(i)))$ for $i=2, \ldots, l$ and
$\textbf{t}''(i) = \textbf{u}''(\pi(i))$ for $i=1, \ldots,l$. The
set of all equivalence classes under the relation $\sim$ is
denoted by $TS = LTS / \sim$. We denote by $\alpha(\textbf{t})$
the cardinality of $\textbf{t}$, i.e.\ the number of possibilities
of monotonically labelling the nodes of $\textbf{t}$ with numbers
$1, \ldots, l(\textbf{t})$.
%\\ \\
%
For example, the labelled trees $[[\tau_{0}^3]_{j_1}^2,
\tau_{0}^4]_{\gamma}^1$, $[[\tau_{0}^4]_{j_1}^2,
\tau_{0}^3]_{\gamma}^1$ and $[\tau_{0}^2,
[\tau_{0}^4]_{j_1}^3]_{\gamma}^1$ with $j_1 \in \{1, \ldots, m\}$
belong to the same equivalence class as $\textbf{t}_I$ in the
example above. Thus, we have $\alpha(\textbf{t}_I)=3$. For $j_1, j_2
\in \{1, \ldots, m\}$ with $j_1 \neq j_2$, we obtain the two
different labelled trees $[[\tau_{j_1}^3, \tau_{j_2}^4]_0^2
]_{\gamma}^1$ and $[[\tau_{j_1}^4, \tau_{j_2}^3]_0^2 ]_{\gamma}^1$
belonging to the same equivalence class as $\textbf{t}_{II}$ and we
get $\alpha([[\tau_{j_1}, \tau_{j_2}]_0 ]_{\gamma})=2$. However, if
we choose $j_1=j_2$, then there exists only on labelled tree
$[[\tau_{j_1}^3, \tau_{j_1}^4]_0^2 ]_{\gamma}^1$ and we obtain
$\alpha([[\tau_{j_1}, \tau_{j_1}]_0 ]_{\gamma})=1$.

%
%\begin{figure}[Htbp]
%\begin{center}
%    \begin{tabular}{ccccc}
%    \pstree[treemode=U, dotstyle=otimes, dotsize=3.2mm, levelsep=0.1cm, radius=1.6mm, treefit=loose]
%    {\Tn}{
%    \pstree[treemode=U, dotstyle=otimes, dotsize=3.2mm, levelsep=1cm, radius=1.6mm, treefit=loose]
%    {\Tdot~[tnpos=l]{1 }} {\pstree{\TC*~[tnpos=l]{2}}{\TC~[tnpos=l]{3}~[tnpos=r]{$\!\!_{j_1}$}}
%    \TC~[tnpos=l]{4}~[tnpos=r]{$\!\!_{j_2}$}}
%    }
%    & \qquad \qquad &
%    \pstree[treemode=U, dotstyle=otimes, dotsize=3.2mm, levelsep=0.1cm, radius=1.6mm, treefit=loose]
%    {\Tn}{
%    \pstree[treemode=U, dotstyle=otimes, dotsize=3.2mm, levelsep=1cm, radius=1.6mm, treefit=loose]
%    {\Tdot~[tnpos=l]{1 }} {\pstree{\TC*~[tnpos=l]{2}}{\TC~[tnpos=l]{4}~[tnpos=r]{$\!\!_{j_2}$}}
%    \TC~[tnpos=l]{3}~[tnpos=r]{$\!\!_{j_1}$}}
%    }
%    & \qquad \qquad &
%    \pstree[treemode=U, dotstyle=otimes, dotsize=3.2mm, levelsep=0.1cm, radius=1.6mm, treefit=loose]
%    {\Tn}{
%    \pstree[treemode=U, dotstyle=otimes, dotsize=3.2mm, levelsep=1cm, radius=1.6mm, treefit=loose]
%    {\Tdot~[tnpos=l]{1 }} {\TC~[tnpos=l]{2}~[tnpos=r]{$\!\!_{j_3}$} \pstree{\TC*~[tnpos=l]{3}}{
%    \TC~[tnpos=l]{4}~[tnpos=r]{$\!\!_{j_8}$}}}
%    }
%    \end{tabular}
%\caption{Trees of the same equivalence class.}
%\label{St-equal-trees:Wm}
%\end{center}
%\end{figure}
%
% =======================================================
%
For every rooted tree $\textbf{t} \in TS$, there exists a
corresponding {\emph{elementary differential}}. The elementary
differential is defined recursively by $F(\gamma)(x) = f(x)$,
$F(\tau_0)(x) = a(x)$ and $F(\tau_j)(x) = b^j(x)$ for single nodes
and by
\begin{equation} \label{St-elementary-differential-F:Wm}
    F(\textbf{t})(x) =
    \begin{cases}
    f^{(k)}(x) \cdot (F(\textbf{t}_1)(x), \ldots, F(\textbf{t}_k)(x)) &
    \text{for } \textbf{t}=[\textbf{t}_1, \ldots, \textbf{t}_k]_{\gamma} \\
    a^{(k)}(x) \cdot (F(\textbf{t}_1)(x), \ldots,
    F(\textbf{t}_k)(x)) & \text{for } \textbf{t}=[\textbf{t}_1, \ldots, \textbf{t}_k]_0 \\
    {b^j}^{(k)}(x) \cdot (F(\textbf{t}_1)(x),
    \ldots, F(\textbf{t}_k)(x)) & \text{for } \textbf{t}=[\textbf{t}_1, \ldots,
    \textbf{t}_k]_j
    \end{cases}
\end{equation}
for a tree $\textbf{t}$ with more than one node and $j \in \{1,
\ldots, m\}$. Here $f^{(k)}$, $a^{(k)}$ and ${b^j}^{(k)}$ define a
symmetric $k$-linear differential operator, and one can choose the
sequence of labelled S-trees $\textbf{t}_1, \ldots, \textbf{t}_k$
in an arbitrary order. For example, the $I$th component of
$a^{(k)} \cdot (F(\textbf{t}_1), \ldots, F(\textbf{t}_k))$ can be
written as
\begin{align*}
    ( a^{(k)} \cdot (F(\textbf{t}_1), \ldots, F(\textbf{t}_k)) )^I
    &= \sum_{J_1, \ldots, J_k=1}^d \frac{\partial^k
    a^I}{\partial x^{J_1} \ldots \partial x^{J_k}} \,
    (F^{J_1}(\textbf{t}_1), \ldots, F^{J_k}(\textbf{t}_k))
\end{align*}
where the components of vectors are denoted by superscript
indices, which are chosen as capitals.
As a result of this we calculate for the trees in
Figure~\ref{St-S-tree-examples-tI+tII:Wm} the elementary
differentials
\begin{equation*}
    \begin{split}
    F(\textbf{t}_I) &= f'' ({b^{j_1}}' (a), a) = \sum_{J_1,J_2=1}^d
    \frac{\partial^2 f}{\partial x^{J_1} \partial x^{J_2}}
    \big( \sum_{K_1=1}^d \frac{\partial b^{J_1,j_1}}{\partial x^{K_1}}
    \, a^{K_1} \cdot a^{J_2} \big) \, , \\
%    , \qquad \qquad
    F(\textbf{t}_{II}) &= f' (a'' (b^{j_1}, b^{j_2})) = \sum_{J_1=1}^d
    \frac{\partial f}{\partial x^{J_1}} \big( \sum_{K_1, K_2 =1}^d
    \frac{\partial^2 a^{J_1}}{\partial x^{K_1} \partial
    x^{K_2}} \, b^{K_1,j_1} \cdot b^{K_2,j_2} \big) \, , \\
    F(\textbf{t}_{III}) &= f'' ({b^{j_1}}'(b^{j_3}),
    {b^{j_2}}'(b^{j_4})) \\
    &= \sum_{J_1,J_2=1}^d
    \frac{\partial^2 f}{\partial x^{J_1} \partial x^{J_2}}
    \big( \sum_{K_1, K_2 =1}^d
    \frac{\partial b^{J_1,j_1}}{\partial x^{K_1}} \, b^{K_1,j_3} \cdot
    \frac{\partial b^{J_2,j_2}}{\partial x^{K_2}} \, b^{K_2,j_4}
    \big) \, .
    \end{split}
\end{equation*}

Next, we assign to every tree a corresponding multiple stochastic
integral. For $\textbf{t} \in TS$ and an adapted right continuous
stochastic process $(Z_t)_{t \geq t_0}$ the corresponding
{\emph{multiple stochastic integral}} is recursively defined by
\begin{equation} \label{multiple-Ito-integral}
    I_{\textbf{t};t_0,t}[Z_{\cdot}] = \begin{cases}
    \displaystyle
    ( \prod_{i=1}^k I_{\textbf{t}_i;t_0,t} ) [Z_{\cdot}] & \text{if }
    \textbf{t} = [\textbf{t}_1, \ldots, \textbf{t}_k]_{\gamma} \\
    \displaystyle
    ( \int_{t_0}^t \prod_{i=1}^k I_{\textbf{t}_i;t_0,s}
    *{\mathrm{d}}W_s^j )[Z_{\cdot}] & \text{if } \textbf{t}=[\textbf{t}_1,
    \ldots, \textbf{t}_k]_j
    \end{cases}
\end{equation}
with $*{\mathrm{d}}W_s^0={\mathrm{d}}s$,
$I_{\tau_j;t_0,t}[Z_{\cdot}] = \int_{t_0}^t Z_s
*{\mathrm{d}}W_s^j$, $I_{\gamma;t_0,t}[Z_{\cdot}] = Z_t$,
$I_{\textbf{t};t_0,t} = I_{\textbf{t};t_0,t}[1]$ and with the
notation
\begin{equation}
    \begin{split}
    &( \int_{t_0}^t \int_{t_0}^{s_n} \cdots \int_{t_0}^{s_2}
    \, *{\mathrm{d}}W_{s_1}^{j_1} \, *{\mathrm{d}}W_{s_2}^{j_2}
    \cdots \, *{\mathrm{d}}W_{s_n}^{j_n} ) [Z_{\cdot}]  = \mathcal{I}_{(j_1, j_2,
    \ldots, j_n)}[Z_{\cdot}]_{t_0,t} \\
    &=
    \int_{t_0}^t \int_{t_0}^{s_n} \cdots \int_{t_0}^{s_2} Z_{s_1}
    \, *{\mathrm{d}}W_{s_1}^{j_1} \, *{\mathrm{d}}W_{s_2}^{j_2}
    \cdots \, *{\mathrm{d}}W_{s_n}^{j_n}
    \end{split}
\end{equation}
in (\ref{multiple-Ito-integral}). The product of two stochastic
integrals can be written as a sum
\begin{equation}
    \begin{split}
    &\int_{t_0}^t X_s \, *{\mathrm{d}}W_s^i \, \int_{t_0}^t Y_s \,
    *{\mathrm{d}}W_s^j = \int_{t_0}^t X_s \, Y_s \, 1_{ \{i = j \neq 0 \wedge * \neq \circ\}}
    \, {\mathrm{d}}s \\
    &+ \int_{t_0}^t X_s \, \int_{t_0}^s Y_u \,
    *{\mathrm{d}}W_u^j \, *{\mathrm{d}}W_s^i + \int_{t_0}^t \int_{t_0}^s
    X_u \, *{\mathrm{d}}W_u^i \, Y_s \, *{\mathrm{d}}W_s^j
    \end{split}
\end{equation}
for $0 \leq i,j \leq m$ \cite{KP99}, where the first summand on the
right hand side appears only in the case of It{\^o} calculus.
%
%On the other hand, for Stratonovich integrals we obtain the rule
%\begin{equation}
%    \begin{split}
%    &\int_{t_0}^t X_s \, \circ~{\mathrm{d}}W_s^i \, \int_{t_0}^t Y_s \,
%    \circ~{\mathrm{d}}W_s^j \\
%    &= \int_{t_0}^t X_s \int_{t_0}^s Y_u \,
%    \circ~{\mathrm{d}}W_u^j \, \circ~{\mathrm{d}}W_s^i + \int_{t_0}^t
%    \int_{t_0}^s X_u \, \circ~{\mathrm{d}}W_u^i \, Y_s \,
%    \circ~{\mathrm{d}}W_s^j
%    \end{split}
%\end{equation}
%for $i,j=0,1, \ldots, m$.
For example, we calculate for $\textbf{t}_I$
\begin{equation}
    \begin{split}
    I_{\textbf{t}_I;t_0,t}[Z_{\cdot}] &= (I_{\tau_0;t_0,t} \,
    I_{[\tau_0]_{j_1};t_0,t})[Z_{\cdot}] = ( \int_{t_0}^t \,
    {\mathrm{d}}s \, \int_{t_0}^t I_{\tau_0;t_0,s} \,
    *{\mathrm{d}}W_s^{j_1} ) [Z_{\cdot}] \\
    &= ( \int_{t_0}^t \int_{t_0}^s I_{\tau_0;t_0,u} \,
    *{\mathrm{d}}W_u^{j_1} \, {\mathrm{d}}s + \int_{t_0}^t \int_{t_0}^s \,
    {\mathrm{d}}u \, I_{\tau_0;t_0,s} \, *{\mathrm{d}}W_s^{j_1} )
    [Z_{\cdot}] \\
    &= ( \int_{t_0}^t \int_{t_0}^s \int_{t_0}^u \, {\mathrm{d}}v \,
    *{\mathrm{d}}W_u^{j_1} \, {\mathrm{d}}s + \int_{t_0}^t \int_{t_0}^s \,
    {\mathrm{d}}u \, \int_{t_0}^s \, {\mathrm{d}}u \, *{\mathrm{d}}W_s^{j_1} )
    [Z_{\cdot}] \\
    &= \mathcal{I}_{(0,j_1,0)}[Z_{\cdot}]_{t_0,t} + 2
    \mathcal{I}_{(0,0,j_1)}[Z_{\cdot}]_{t_0,t} \, .
%    &= \int_{t_0}^t \int_{t_0}^s \int_{t_0}^u Z_{v} \, {\mathrm{d}}v \,
%    *{\mathrm{d}}W_u^{j_1} \, {\mathrm{d}}s + 2 \int_{t_0}^t \int_{t_0}^s
%    \, \int_{t_0}^u Z_{v} \, {\mathrm{d}}v \, {\mathrm{d}}u
%    \, *{\mathrm{d}}W_s^{j_1} \, .
    \end{split}
\end{equation}
For the tree $\textbf{t}_{II}$ we obtain
\begin{equation}
    \begin{split}
    I_{\textbf{t}_{II};t_0,t}[Z_{\cdot}] &= (\int_{t_0}^t
    I_{\tau_{j_1};t_0,s} \, I_{\tau_{j_2};t_0,s} \, {\mathrm{d}}s)
    [Z_{\cdot}] \\
    &= (\int_{t_0}^t \int_{t_0}^s
    I_{\tau_{j_2};t_0,u} \, *{\mathrm{d}}W_u^{j_1} \, {\mathrm{d}}s
    + \int_{t_0}^t \int_{t_0}^s
    I_{\tau_{j_1};t_0,u} \, *{\mathrm{d}}W_u^{j_2} \, {\mathrm{d}}s
    \\
    &+ \int_{t_0}^t \int_{t_0}^s
    1_{\{j_1=j_2 \neq 0 \wedge * \neq \circ\}} \, {\mathrm{d}}u \, {\mathrm{d}}s)
    [Z_{\cdot}] \\
    &= \mathcal{I}_{(j_2,j_1,0)}[Z_{\cdot}]_{t_0,t} +
    \mathcal{I}_{(j_1,j_2,0)}[Z_{\cdot}]_{t_0,t} +
    \mathcal{I}_{(0,0)}[Z_{\cdot} \, 1_{\{j_1=j_2 \neq 0 \wedge * \neq
    \circ\}}]_{t_0,t} \, .
%    &= \int_{t_0}^t \int_{t_0}^s
%    \int_{t_0}^u Z_v \, *{\mathrm{d}}W_v^{j_2} \, *{\mathrm{d}}W_u^{j_1} \, {\mathrm{d}}s
%    + \int_{t_0}^t \int_{t_0}^s
%    \int_{t_0}^u Z_v \, *{\mathrm{d}}W_v^{j_1} \, *{\mathrm{d}}W_u^{j_2} \, {\mathrm{d}}s
%    \\
%    &+ \int_{t_0}^t \int_{t_0}^s
%    Z_u \, 1_{\{j_1=j_2 \neq 0 \wedge * \neq \circ\}} \, {\mathrm{d}}u \,
%    {\mathrm{d}}s \, .
    \end{split}
\end{equation}
%Analogously, we obtain for $\textbf{t}_{III}$ the corresponding
%multiple integral
%\begin{equation}
%    \begin{split}
%    I_{\textbf{t}_{III};t_0,t}[Z_{\cdot}] &= (\int_{t_0}^t \int_{t_0}^s
%    \, *{\mathrm{d}}W_u^{j_3} \, *{\mathrm{d}}W_s^{j_1} \,
%    \int_{t_0}^t \int_{t_0}^s \, *{\mathrm{d}}W_u^{j_4} \,
%    *{\mathrm{d}}W_s^{j_2})[Z_{\cdot}] \\
%    &= \sum_{\sigma} \int_{t_0}^t \int_{t_0}^{s_4} \int_{t_0}^{s_3}
%    \int_{t_0}^{s_2} Z_{s_1} \, *{\mathrm{d}}W_{s_1}^{\sigma(j_1)} \,
%    *{\mathrm{d}}W_{s_2}^{\sigma(j_2)}
%    \, *{\mathrm{d}}W_{s_3}^{\sigma(j_3)} \,
%    *{\mathrm{d}}W_{s_4}^{\sigma(j_4)} + ...
%    \end{split}
%\end{equation}
%where the sum is taken over all permutations $\sigma$ of
%$(j_1,j_2,j_3,j_4)$ such that .... .

Let $\textbf{t} \in TS$ with $\textbf{t} = [\textbf{t}_1, \ldots,
\textbf{t}_1, \textbf{t}_2, \ldots , \textbf{t}_2, \ldots,
\textbf{t}_k, \ldots, \textbf{t}_k]_j = [\textbf{t}_1^{n_1},
\textbf{t}_2^{n_2}, \ldots, \textbf{t}_k^{n_k}]_j$, $j \in \{\gamma,
0, 1, \ldots, m\}$,
%or $\textbf{t} = (\textbf{t}_1, \ldots, \textbf{t}_1,
%\textbf{t}_2, \ldots , \textbf{t}_2, \ldots, \textbf{t}_k, \ldots,
%\textbf{t}_k) = (\textbf{t}_1^{n_1}, \textbf{t}_2^{n_2}, \ldots,
%\textbf{t}_k^{n_k})$
where $\textbf{t}_1, \ldots, \textbf{t}_k$ are distinct subtrees
with multiplicities $n_1, \ldots, n_k$, respectively. Then we
recursively define the {\em symmetry factor} by
\begin{equation} \label{sym-factor}
    \sigma(\textbf{t}) = \prod_{i=1}^k n_i! \,\,
    \sigma(\textbf{t}_i)^{n_i} \, .
\end{equation}
Next, we define the density of a tree which is a measure of its
non-bushiness. For $\textbf{t} \in TS$ let the density
$\gamma(\textbf{t})$ be recursively defined by
$\gamma(\textbf{t})=1$ if $l(\textbf{t})=1$ and
\begin{equation}
    \gamma(\textbf{t}) =
    l(\textbf{t}) \prod_{i=1}^{k} \gamma(\textbf{t}_i)
\end{equation}
if $\textbf{t}=[\textbf{t}_1, \ldots, \textbf{t}_{k}]_j$ with some
$j \in \{\gamma, 0,1, \ldots, m\}$. Then, for $\textbf{t} \in TS$ we
obtain
\begin{equation}
    \alpha(\textbf{t}) = \frac{l(\textbf{t})!}{\gamma(\textbf{t}) \,
    \sigma(\textbf{t})}
\end{equation}
see also \cite{Butcher03}. For example, we calculate for
$\textbf{t}_I$ that $\sigma(\textbf{t}_I)=1$ and with
$l(\textbf{t}_I)!=24$ and $\gamma(\textbf{t}_I)=8$ we obtain
$\alpha(\textbf{t}_I)=3$. For the tree $\textbf{t}_{II}$ we have
to consider two cases: if $j_1 \neq j_2$ we have
$\sigma(\textbf{t}_{II})=1$ and we get with
$l(\textbf{t}_{II})!=24$ and $\gamma(\textbf{t}_{II})=12$ that
$\alpha(\textbf{t}_{II})=2$. However, in the case of $j_1=j_2$ we
have some symmetry and thus $\sigma(\textbf{t}_{II})=2$ which
results in $\alpha(\textbf{t}_{II})=1$. Analogously, we consider
for $\textbf{t}_{III}$ the case of $j_1=j_2$ and $j_3=j_4$ with
some symmetry where $\sigma(\textbf{t}_{III})=2$ and with
$l(\textbf{t})!=120$ and $\gamma(\textbf{t}_{III})=20$ follows
$\alpha(\textbf{t}_{III})=3$. For all other cases, we have no
symmetry and thus $\sigma(\textbf{t}_{III})=1$ where we obtain
$\alpha(\textbf{t}_{III})=6$ different monotonically labelled
trees in the equivalence class of $\textbf{t}_{III}$.
%
%
%
%
% ============================================================================
\section{Stochastic Taylor Expansion}
\label{Sec:Taylor-Expansion-SDE}
% ============================================================================
%
%
In order to give a stochastic Taylor expansion of the solution of
the considered SDE~(\ref{Intro-Ito-St-SDE1-integralform-Wm}) with
some remainder term, we have to introduce the sets of descendant
trees. For $\textbf{t} \in TS$ and $j \in \{0,1, \ldots, m\}$ let
$H^j(\textbf{t})$ denote the set of all trees in $TS$ which are
obtained from $\textbf{t}$ by adding one node $\tau_j$. Further,
let $H^{I}(\textbf{t})$ denote the set of all trees in $TS$ which
are obtained from $\textbf{t}$ by adding the two nodes $\tau_j^a$
and $\tau_j^b$ where both nodes have the same $j \in \{1, \ldots,
m\}$, $a=l(\textbf{t})+1$, $b=l(\textbf{t})+2$ and where neither
of them is father of the other. For example, if $\textbf{t} =
[\tau_{j_1},\tau_{j_2}]_{\gamma}$ for some arbitrarily fixed $j_1,
j_2 \in \{1, \ldots, m\}$ then we obtain
\begin{equation}
    \begin{split}
    H^0(\textbf{t}) = \{ &[\tau_{j_1}, \tau_{j_2}, \tau_0]_{\gamma},
    [\tau_{j_1},[\tau_0]_{j_2}]_{\gamma},
    [[\tau_0]_{j_1},\tau_{j_2}]_{\gamma} \} \, , \\
    H^I(\textbf{t}) = \{ &[\tau_{j_1}, \tau_{j_2}, \tau_{j_3}^a,
    \tau_{j_3}^b]_{\gamma},
    [[\tau_{j_3}^a, \tau_{j_3}^b]_{j_1}, \tau_{j_2}]_{\gamma},
    [\tau_{j_1},[\tau_{j_3}^a, \tau_{j_3}^b]_{j_2}]_{\gamma},
    [[\tau_{j_3}^a]_{j_1}, \tau_{j_2}, \tau_{j_3}^b]_{\gamma}, \\
    &[[\tau_{j_3}^b]_{j_1}, \tau_{j_2}, \tau_{j_3}^a]_{\gamma},
    [\tau_{j_1}, [\tau_{j_3}^a]_{j_2}, \tau_{j_3}^b]_{\gamma},
    [\tau_{j_1}, [\tau_{j_3}^b]_{j_2}, \tau_{j_3}^a]_{\gamma}, \\
    &[[\tau_{j_3}^a]_{j_1}, [\tau_{j_3}^b]_{j_2}]_{\gamma},
    [[\tau_{j_3}^b]_{j_1}, [\tau_{j_3}^a]_{j_2}]_{\gamma}
    : j_3^a=j_3^b \in \{1, \ldots, m\} \}
    \end{split}
\end{equation}
independently whether $j_1=j_2$ or $j_1 \neq j_2$. In the
following, let $\tfrac{1}{2} \mathbb{N}_0 = \{p: 2p \in
\mathbb{N}_0\}$.
Then, based on the introduced multi-colored rooted trees, we
obtain the following stochastic Taylor expansion for the solution
of the It{\^o} SDE~(\ref{Intro-Ito-St-SDE1-integralform-Wm}):

\begin{The} \label{Ito-tree-expansion-exact-sol:Wm}
    For the solution process $(X_t)_{t \in [t_0,T]}$ of the It{\^o}
    SDE~(\ref{Intro-Ito-St-SDE1-integralform-Wm}) and for
    $p \in \tfrac{1}{2} \mathbb{N}_0$, $f :
    \mathbb{R}^d \rightarrow \mathbb{R}$
    with $f, a^i, b^{i,j}
    \in C^{2p+2}(\mathbb{R}^d, \mathbb{R})$
    for $i=1, \ldots, d$, $j=1, \ldots, m$, we obtain the
    expansion
    \begin{equation}
    \begin{split} \label{Ito-Tree-expansion-exact-sol-formula1:Wm}
        f(X_{t})
        = &\sum_{\substack{\textbf{t} \in TS \\ \rho(\textbf{t}) \leq
        p}} F(\textbf{t})(X_{t_0}) \, \frac{I_{\textbf{t};t_0,t}}{\sigma(\textbf{t})}
        + \mathcal{R}_{p}(t,t_0)
    \end{split}
    \end{equation}
    $\Prob$-a.s.\ with remainder term
    \begin{equation}
    \label{Ito-Tree-expansion-exact-sol-remainder1:Wm}
        \begin{split}
        \mathcal{R}_{p}(t,t_0) &=
        \sum_{\substack{\textbf{t} \in TS \\ \rho(\textbf{t}) = p+1/2}}
        \frac{I_{\textbf{t};t_0,t}[F(\textbf{t})(X_{\cdot})]}{\sigma(\textbf{t})}
        + \sum_{\substack{\textbf{t} \in TS \\ \rho(\textbf{t}) = p}}
        \sum_{\textbf{u} \in H^0(\textbf{t})}
        \frac{I_{\textbf{t};t_0,t}[\int_{t_0}^{\cdot}
        F(\textbf{u})(X_{s})
        \, {\mathrm{d}}s]}{\sigma(\textbf{t})} \\
        &\quad + \sum_{\substack{\textbf{t} \in TS \\ \rho(\textbf{t}) = p}}
        \sum_{\textbf{u} \in H^I(\textbf{t})}
        \frac{I_{\textbf{t};t_0,t}[ \int_{t_0}^{\cdot}
        F(\textbf{u})(X_{s})
        \, {\mathrm{d}}s]}{2 \, \sigma(\textbf{t})}
        \end{split}
    \end{equation}
    provided all of the appearing multiple It{\^o} integrals exist.
\end{The}

{\bf{Proof.}} First, we assign to every $\textbf{t} \in LTS$ a
corresponding multiple stochastic integral. Therefore, let $\Gamma
= \{ \textbf{t} \in LTS : \textbf{t}''(l(\textbf{t})) = \tau_0
\}$, $\Lambda^j = \{ \textbf{t} \in LTS :
\textbf{t}''(l(\textbf{t})) = \textbf{t}''(l(\textbf{t})-1) =
\tau_j \, \wedge \, \textbf{t}'(l(\textbf{t})) \neq
l(\textbf{t})-1\}$ and $\Sigma^j = \{ \textbf{t} \in LTS : (
\textbf{t}''(l(\textbf{t})) = \tau_j \, \wedge \,
\textbf{t}''(l(\textbf{t})-1) \neq \tau_j ) \, \vee \, (
\textbf{t}''(l(\textbf{t})) = \textbf{t}''(l(\textbf{t})-1) =
\tau_j \, \wedge \, \textbf{t}'(l(\textbf{t})) = l(\textbf{t})-1
\}$ for $j=1, \ldots, m$. Then, define for $\textbf{t} \in LTS$
with $l=l(\textbf{t})$ nodes and for an adapted right continuous
stochastic process $(Z_t)_{t \geq t_0}$ the corresponding multiple
It{\^o} integral recursively by
\begin{equation}
    \hat{I}_{\textbf{t}}[Z_{\cdot}]_{t_0,t} = \begin{cases} Z_t & \text{if }
    \rho(t)=0 \\
    \hat{I}_{\textbf{t}_{-}} [\int_{t_0}^{\cdot} Z_s \, {\mathrm{d}}s ]_{t_0,t} & \text{if }
    \textbf{t} \in \Gamma \\
    \hat{I}_{\textbf{t}_{-}} [\int_{t_0}^{\cdot} Z_s \,
    {\mathrm{d}}W^{j}_s
    ]_{t_0,t} + \hat{I}_{\textbf{t}_{--}} [\frac{1}{2} \int_{t_0}^{\cdot}
    Z_s \, {\mathrm{d}}s ]_{t_0,t}
    & \text{if } \textbf{t} \in \Lambda^j \\
    \hat{I}_{\textbf{t}_{-}} [\int_{t_0}^{\cdot} Z_s \,
    {\mathrm{d}}W^{j}_s ]_{t_0,t} & \text{if } \textbf{t} \in \Sigma^j
    \end{cases}
\end{equation}
where $\textbf{t}_{-}$ is the tree which is obtained from
$\textbf{t}$ by removing the last node with label $l(\textbf{t})$
and $\textbf{t}_{--} = (\textbf{t}_{-})_{-}$ denotes the tree
where the last two nodes with labels $l(\textbf{t})$ and
$l(\textbf{t})-1$ are removed.
%
%For $j=1, \ldots, m$, let
%\begin{equation}
%    \hat{L}^{0} = \sum_{k=1}^d a^k \frac{\partial}{\partial x^k},
%    \quad \quad
%    \hat{L}^{j} = \sum_{k,l=1}^d b^{k,j} \, b^{l,j}
%    \frac{\partial^2}{\partial x^k \partial x^l}, \quad \quad
%    L^j = \sum_{k=1}^d b^{k,j} \frac{\partial}{\partial x^k},
%\end{equation}
%and $L^0 = \hat{L}^{0} + \tfrac{1}{2} \sum_{j=1}^m \hat{L}^{j}$.
%
Following the notation in \cite{KP99}, for a multi-index $\alpha =
(j_1, \ldots, j_l) \in \{0,1, \ldots, m\}^l$ let $l(\alpha)=l$ be
the length with $l(\nu)=0$ for the multi-index $\nu$ of length 0.
Further, let $\mathcal{M}$ be the set of all multi-indices and
$n(\alpha)$ be the number of components of $\alpha$ which are
equal to $0$. For $j \in \{0,1, \ldots, m\}$ let $(j)
* \alpha = (j, j_1, \ldots, j_l)$. Then, define $\mathcal{I}_{(j_1,
\ldots, j_l)}[Z_s]_{t_0,t} = \int_{t_0}^t \mathcal{I}_{(j_1,
\ldots, j_{l-1})}[Z_u]_{t_0,s} \, {\mathrm{d}}W^{j_l}_s$ if $l
\geq 1$ and $\mathcal{I}_{\nu}[Z_s]_{t_0,t} = Z_t$. Finally, let
$f_{(j_1, \ldots, j_l)} = L^{j_1} f_{(j_2, \ldots, j_l)}$ and
$f_{\nu}=f$.
In the following, we consider the hierarchical set $\mathcal{A}_p
= \{ \alpha \in \mathcal{M} : l(\alpha)+n(\alpha) \leq 2p\}$ and
define $A_p = \{\alpha \in \mathcal{M} : l(\alpha)+n(\alpha)=2p\}$
for $p \in \frac{1}{2} \mathbb{N}_0$ with $\mathcal{A}_p =
\bigcup_{i \in \frac{1}{2} \mathbb{N}_0 : i \leq p} A_i$ and $A_i
\cap A_j = \emptyset$ if $i \neq j$. Then, a recursive application
of the It{\^o}-formula yields the It{\^o}-Taylor expansion~\cite{KP99}
\begin{equation} \label{Proof-Ito-Taylor-Formel}
    f(X_t) = \sum_{\alpha \in \mathcal{A}_p}
    \mathcal{I}_{\alpha}[f_{\alpha}(X_{t_0})]_{t_0,t} +
    \sum_{\alpha \in \mathcal{B}(\mathcal{A}_{p})}
    \mathcal{I}_{\alpha}[f_{\alpha}(X_{\cdot})]_{t_0,t}
\end{equation}
$\Prob$-a.s.\ with $\mathcal{B}(\mathcal{A}_p) = \{ \alpha \in
\mathcal{M} \setminus \mathcal{A}_p : \alpha = (j) * \bar{\alpha},
j \in \{0,1,\ldots,m\}, \bar{\alpha} \in \mathcal{A}_p \}$.
Clearly, $A_p = \{(0)*\alpha \in \mathcal{M} : \alpha \in
A_{p-1}\} \cup \{(j)*\alpha \in \mathcal{M} : \alpha \in
A_{p-\frac{1}{2}}, j \in \{1, \ldots, m\}\}$ for $p \geq 1$. As a
result of this, we have to prove
\begin{equation} \label{Proof-zz-Hauptteil}
    \sum_{\alpha \in A_{\tfrac{1}{2}n}}
    \mathcal{I}_{\alpha}[f_{\alpha}(X_{t_0})]_{t_0,t}
    = \sum_{\substack{\textbf{t} \in LTS \\ \rho(\textbf{t}) = \tfrac{1}{2} n}}
    F(\textbf{t})(X_{t_0}) \, \hat{I}_{\textbf{t}}[1]_{t_0,t}
\end{equation}
for $n \in \mathbb{N}_0$. Now, define a linear operator $K^j$ for
$j=0,1, \ldots, m$ by
\begin{equation} \label{Proof-linearer-Operator}
    K^j(\mathcal{I}_{(j_1, \ldots, j_l)}[Z_{\cdot}]_{t_0,t})
    = \mathcal{I}_{(j_1, \ldots, j_l)}[\int_{t_0}^{\cdot} L^j Z_s \, {\mathrm{d}}W^{j}_s]_{t_0,t}
    = \mathcal{I}_{(j,j_1, \ldots, j_l)}[L^j \, Z_{\cdot}]_{t_0,t} \, .
\end{equation}
Then, for $n=0$, we obtain
\begin{equation}
    \begin{split}
    \sum_{\alpha \in A_{0}}
    \mathcal{I}_{\alpha}[f_{\alpha}(X_{t_0})]_{t_0,t}
    &= \mathcal{I}_{\nu}[f_{\nu}(X_{t_0})]_{t_0,t}
    = f(X_{t_0})
    = F(\gamma)(X_{t_0}) \, \hat{I}_{\gamma}[1]_{t_0,t} \\
    &= \sum_{\substack{\textbf{t} \in LTS \\ \rho(\textbf{t}) = 0}}
    F(\textbf{t})(X_{t_0}) \hat{I}_{\textbf{t}}[1]_{t_0,t}
    \end{split}
\end{equation}
and for $n=1$, we get
\begin{equation}
    \begin{split}
    \sum_{\alpha \in A_{\tfrac{1}{2}}}
    \mathcal{I}_{\alpha}[f_{\alpha}(X_{t_0})]_{t_0,t}
    &= \sum_{j=1}^m \mathcal{I}_{(j)}[f_{(j)}(X_{t_0})]_{t_0,t} \\
    &= \sum_{j=1}^m \sum_{k=1}^d \frac{\partial f}{\partial
    x^k}(X_{t_0}) \, b^{k,j}(X_{t_0}) \,
    \int_{t_0}^t 1 \, {\mathrm{d}}W_s^j \\
    &= \sum_{j=1}^m F([\tau_j]_{\gamma})(X_{t_0}) \,
    \hat{I}_{[\tau_j]_{\gamma}}[1]_{t_0,t} \\
    &= \sum_{\substack{\textbf{t} \in LTS \\ \rho(\textbf{t}) =
    \frac{1}{2}}} F(\textbf{t})(X_{t_0}) \hat{I}_{\textbf{t}}[1]_{t_0,t}
    \end{split}
\end{equation}
Now, assume that (\ref{Proof-zz-Hauptteil}) holds for some $n,n-1
\geq 0 $. For $\textbf{t} \in LTS$ and $j=0,1, \ldots, m$, we
introduce the sets
\begin{equation}
    \begin{split}
    \hat{H}_1^j(\textbf{t}) = \{ & \textbf{u} \in LTS :
    l(\textbf{u})=l(\textbf{t})+1, \textbf{u}' |_{\{2, \ldots,
    l(\textbf{t})\}}=\textbf{t}', \textbf{u}'' |_{\{1, \ldots,
    l(\textbf{t})\}}=\textbf{t}'', \\
    & \textbf{u}''(l(\textbf{t})+1)=\tau_j,
    \textbf{u}'(l(\textbf{t})+1) \neq l(\textbf{t}) \} \, , \\
    \bar{H}_1^j(\textbf{t}) = \{ & \textbf{u} \in LTS :
    l(\textbf{u})=l(\textbf{t})+1, \textbf{u}' |_{\{2, \ldots,
    l(\textbf{t})\}}=\textbf{t}', \textbf{u}'' |_{\{1, \ldots,
    l(\textbf{t})\}}=\textbf{t}'', \\
    & \textbf{u}''(l(\textbf{t})+1)=\tau_j,
    \textbf{u}'(l(\textbf{t})+1) = l(\textbf{t}) \} \, , \\
    \hat{H}_2^j(\textbf{t}) = \{ & \textbf{u} \in LTS :
    l(\textbf{u})=l(\textbf{t})+2, \textbf{u}' |_{\{2, \ldots,
    l(\textbf{t})\}}=\textbf{t}', \textbf{u}'' |_{\{1, \ldots,
    l(\textbf{t})\}}=\textbf{t}'', \\
    & \textbf{u}''(l(\textbf{t})+1)=\textbf{u}''(l(\textbf{t})+2)=\tau_j,
    \textbf{u}'(l(\textbf{t})+2) \neq l(\textbf{t})+1 \} \, , \\
    \bar{H}_2^j(\textbf{t}) = \{ & \textbf{u} \in LTS :
    l(\textbf{u})=l(\textbf{t})+2, \textbf{u}' |_{\{2, \ldots,
    l(\textbf{t})\}}=\textbf{t}', \textbf{u}'' |_{\{1, \ldots,
    l(\textbf{t})\}}=\textbf{t}'', \\
    & \textbf{u}''(l(\textbf{t})+1)=\textbf{u}''(l(\textbf{t})+2)=\tau_j,
    \textbf{u}'(l(\textbf{t})+2) = l(\textbf{t})+1 \} \, ,
    \end{split}
\end{equation}
and let $H_i^j = \hat{H}_i^j \cup \bar{H}_i^j$ for $i=1,2$.
Then, with Lemma~2.7 and Lemma~2.8 in \cite{Roe04b} we obtain
{\allowdisplaybreaks
%\begin{equation*}
%    \begin{split}
\begin{alignat*}{3}
    &\sum_{\alpha \in A_{{1}/{2}(n+1)}}
    \mathcal{I}_{\alpha}[f_{\alpha}(X_{t_0})]_{t_0,t} \\
    &= K^0 (\sum_{\alpha \in A_{{1}/{2}(n-1)}}
    \mathcal{I}_{\alpha}[f_{\alpha}(X_{t_0})]_{t_0,t}) +
    \sum_{j=1}^m K^j (\sum_{\alpha \in A_{{1}/{2}n}}
    \mathcal{I}_{\alpha}[f_{\alpha}(X_{t_0})]_{t_0,t}) \\
    &= K^0 (\sum_{\substack{\textbf{t} \in LTS \\ \rho(\textbf{t}) = {1}/{2}(n-1)}}
    \hat{I}_{\textbf{t}}[F(\textbf{t})(X_{t_0})]_{t_0,t}) +
    \sum_{j=1}^m K^j (\sum_{\substack{\textbf{t} \in LTS \\ \rho(\textbf{t}) = {1}/{2}n}}
    \hat{I}_{\textbf{t}}[F(\textbf{t})(X_{t_0})]_{t_0,t}) \\
%    &= \sum_{\substack{\textbf{t} \in LTS \\ \rho(\textbf{t}) = {1}/{2}(n-1)}}
%    \hat{I}_{\textbf{t}}[\int_{t_0}^{\cdot} L^0 F(\textbf{t})(X_{t_0}) \, {\mathrm{d}}s]_{t_0,t} \\
%    &\quad +
%    \sum_{j=1}^m \sum_{\substack{\textbf{t} \in LTS \\ \rho(\textbf{t}) = {1}/{2}n}}
%    \hat{I}_{\textbf{t}}[\int_{t_0}^{\cdot} L^j F(\textbf{t})(X_{t_0}) \, {\mathrm{d}}W^{j}_s]_{t_0,t} \\
%
    &= \sum_{\substack{\textbf{t} \in LTS \\ \rho(\textbf{t}) = {1}/{2}(n-1)}}
    \hat{I}_{\textbf{t}}[\int_{t_0}^{\cdot} \hat{L}^0 F(\textbf{t})(X_{t_0}) \,
    {\mathrm{d}}s]_{t_0,t} \\
    &\quad + \sum_{\substack{\textbf{t} \in LTS \\ \rho(\textbf{t}) =
    {1}/{2}(n-1)}} \sum_{j=1}^m
    \hat{I}_{\textbf{t}}[\frac{1}{2} \int_{t_0}^{\cdot} \hat{L}^j F(\textbf{t})(X_{t_0}) \,
    {\mathrm{d}}s]_{t_0,t} \\
    &\quad + \sum_{j=1}^m \sum_{\substack{\textbf{t} \in LTS \\ \rho(\textbf{t}) = {1}/{2}n \\
    \textbf{t}''(l(\textbf{t})) = \tau_j}}
    \hat{I}_{\textbf{t}}[\int_{t_0}^{\cdot} L^j F(\textbf{t})(X_{t_0}) \,
    {\mathrm{d}}W^{j}_s]_{t_0,t} \\
    &\quad + \sum_{j=1}^m \sum_{\substack{\textbf{t} \in LTS \\ \rho(\textbf{t}) = {1}/{2}n \\
    \textbf{t}''(l(\textbf{t})) \neq \tau_j}}
    \hat{I}_{\textbf{t}}[\int_{t_0}^{\cdot} L^j F(\textbf{t})(X_{t_0}) \, {\mathrm{d}}W^{j}_s]_{t_0,t} \\
    &= \sum_{\substack{\textbf{t} \in LTS \\ \rho(\textbf{t}) =
    {1}/{2}(n-1)}} \sum_{\textbf{u} \in H_1^0(\textbf{t})}
    \hat{I}_{\textbf{u}}[F(\textbf{u})(X_{t_0})]_{t_0,t} \\
    &\quad + \sum_{\substack{\textbf{t} \in LTS \\ \rho(\textbf{t}) =
    {1}/{2}(n-1)}} \sum_{j=1}^m \sum_{\textbf{u} \in
    \hat{H}_2^j(\textbf{t})}
    \hat{I}_{\textbf{t}}[\frac{1}{2} \int_{t_0}^{\cdot} F(\textbf{u})(X_{t_0}) \,
    {\mathrm{d}}s]_{t_0,t} \\
    &\quad + \sum_{j=1}^m \sum_{\substack{\textbf{t} \in LTS \\ \rho(\textbf{t}) = {1}/{2}n \\
    \textbf{t}''(l(\textbf{t})) = \tau_j}}
    ( \sum_{\textbf{u} \in \hat{H}_1^j(\textbf{t})} \hat{I}_{\textbf{t}}[\int_{t_0}^{\cdot} F(\textbf{u})(X_{t_0}) \,
    {\mathrm{d}}W^{j}_s]_{t_0,t} \\
    &\quad + \sum_{\textbf{u} \in
    \bar{H}_1^j(\textbf{t})} \hat{I}_{\textbf{t}}[\int_{t_0}^{\cdot} F(\textbf{u})(X_{t_0}) \,
    {\mathrm{d}}W^{j}_s]_{t_0,t} ) \\
    &\quad + \sum_{j=1}^m \sum_{\substack{\textbf{t} \in LTS \\ \rho(\textbf{t}) = {1}/{2}n \\
    \textbf{t}''(l(\textbf{t})) \neq \tau_j}} \sum_{\textbf{u} \in
    \hat{H}_1^j(\textbf{t}) \cup \bar{H}_1^j(\textbf{t})}
    \hat{I}_{\textbf{u}}[F(\textbf{u})(X_{t_0})]_{t_0,t} \\
    &= \sum_{\substack{\textbf{t} \in LTS \\ \rho(\textbf{t}) =
    {1}/{2}(n-1)}} \sum_{\textbf{u} \in H_1^0(\textbf{t})}
    \hat{I}_{\textbf{u}}[F(\textbf{u})(X_{t_0})]_{t_0,t} \\
    &\quad + \sum_{j=1}^m \sum_{\substack{\textbf{t} \in LTS \\ \rho(\textbf{t}) = {1}/{2}n \\
    \textbf{t}''(l(\textbf{t})) = \tau_j}}
    ( \sum_{\textbf{u} \in \hat{H}_1^j(\textbf{t})} \hat{I}_{\textbf{u}}[F(\textbf{u})(X_{t_0})]_{t_0,t}
    + \sum_{\textbf{u} \in \bar{H}_1^j(\textbf{t})}
    \hat{I}_{\textbf{u}}[F(\textbf{u})(X_{t_0})]_{t_0,t} ) \\
    &\quad + \sum_{j=1}^m \sum_{\substack{\textbf{t} \in LTS \\ \rho(\textbf{t}) = {1}/{2}n \\
    \textbf{t}''(l(\textbf{t})) \neq \tau_j}} \sum_{\textbf{u} \in
    \hat{H}_1^j(\textbf{t}) \cup \bar{H}_1^j(\textbf{t})}
    \hat{I}_{\textbf{u}}[F(\textbf{u})(X_{t_0})]_{t_0,t} \\
    &= \sum_{\substack{\textbf{t} \in LTS \\ \rho(\textbf{t}) =
    {1}/{2}(n-1)}} \sum_{\textbf{u} \in H_1^0(\textbf{t})}
    \hat{I}_{\textbf{u}}[F(\textbf{u})(X_{t_0})]_{t_0,t} \\
    &\quad + \sum_{j=1}^m \sum_{\substack{\textbf{t} \in LTS \\ \rho(\textbf{t}) = {1}/{2}n}}
    \sum_{\textbf{u} \in \hat{H}_1^j(\textbf{t}) \cup \bar{H}_1^j(\textbf{t})}
    \hat{I}_{\textbf{u}}[F(\textbf{u})(X_{t_0})]_{t_0,t} \\
    &= \sum_{\substack{\textbf{t} \in LTS \\ \rho(\textbf{t}) = {1}/{2}(n+1)}}
    \hat{I}_{\textbf{t}}[F(\textbf{t})(X_{t_0})]_{t_0,t}
\end{alignat*}
%    \end{split}
%\end{equation}
}
because $\{\textbf{t} \in LTS : \rho(\textbf{t})=1/2(n+1) \} =
\bigcup_{j=1}^m \{ \textbf{u} \in H_1^j(\textbf{t}) \cup
\hat{H}_1^j(\textbf{t}) : \textbf{t} \in LTS, \rho(\textbf{t}) =
1/2n \} \cup \{ \textbf{u} \in H_1^0(\textbf{t}) : \textbf{t} \in
LTS, \rho(\textbf{t}) = 1/2(n-1) \}$. Thus,
(\ref{Proof-zz-Hauptteil}) holds for all $n \in \mathbb{N}_0$.
%\\ \\

Due to $\mathcal{B}(\mathcal{A}_p) = \bigcup_{j=0}^m \{ (j)*\alpha
\in \mathcal{M} : \alpha \in A_p\} \cup \{ (0)*\alpha \in
\mathcal{M} : \alpha \in A_{p-1/2} \}$ follows analogously for the
remainder term with $p \in \frac{1}{2} \mathbb{N}$
{\allowdisplaybreaks
%\begin{equation}
%    \begin{split}
\begin{alignat*}{3}
    &\sum_{\alpha \in \mathcal{B}(\mathcal{A}_{p})}
    \mathcal{I}_{\alpha}[f_{\alpha}(X_{\cdot})]_{t_0,t}
    = \sum_{\substack{\textbf{t} \in LTS \\ \rho(\textbf{t}) = p}}
    \hat{I}_{\textbf{t}}[\int_{t_0}^{\cdot} L^0 \, F(\textbf{t})(X_{s}) \,
    {\mathrm{d}}s]_{t_0,t} \\
    &\quad + \sum_{\substack{\textbf{t} \in LTS \\ \rho(\textbf{t}) = p}}
    \sum_{j=1}^m \sum_{\textbf{u} \in H_1^j(\textbf{t})}
    \hat{I}_{\textbf{u}}[F(\textbf{u})(X_{\cdot})]_{t_0,t}
%    \\
%    &+ \sum_{\substack{\textbf{t} \in LTS \\ \rho(\textbf{t}) = p}}
%    \sum_{j=1}^m \sum_{\textbf{u} \in \hat{H}_1^j(\textbf{t}) \cup
%    \bar{H}_1^j(\textbf{t})}
%    \hat{I}_{\textbf{u}}[F(\textbf{u})(X_{s})]_{t_0,t}
    + \sum_{\substack{\textbf{t} \in LTS \\ \rho(\textbf{t}) =
    p-1/2}} \sum_{\textbf{u} \in H_1^0(\textbf{t})}
    \hat{I}_{\textbf{u}}[F(\textbf{u})(X_{\cdot})]_{t_0,t} \\
    &= \sum_{\substack{\textbf{t} \in LTS \\ \rho(\textbf{t}) = p}}
%    \sum_{\textbf{u} \in H_1^0(\textbf{t})}
    \hat{I}_{\textbf{t}}[\int_{t_0}^{\cdot} (\hat{L}^0 + \frac{1}{2} \sum_{j=1}^m
    \hat{L}^j) \, F(\textbf{t})(X_{s}) \,
    {\mathrm{d}}s]_{t_0,t}
    + \sum_{\substack{\textbf{t} \in LTS \\ \rho(\textbf{t}) = p+1/2}}
    \hat{I}_{\textbf{t}}[F(\textbf{t})(X_{\cdot})]_{t_0,t}
    \\
    &= \sum_{\substack{\textbf{t} \in LTS \\ \rho(\textbf{t}) = p}}
    \sum_{\textbf{u} \in H_1^0(\textbf{t})}
    \hat{I}_{\textbf{t}}[\int_{t_0}^{\cdot} F(\textbf{u})(X_{s}) \,
    {\mathrm{d}}s]_{t_0,t}
    + \sum_{\substack{\textbf{t} \in LTS \\ \rho(\textbf{t}) = p}}
    \sum_{j=1}^m \sum_{\textbf{u} \in \hat{H}_2^j(\textbf{t})}
    \hat{I}_{\textbf{t}}[\frac{1}{2} \int_{t_0}^{\cdot} F(\textbf{u})(X_{s}) \,
    {\mathrm{d}}s]_{t_0,t} \\
    &+ \sum_{\substack{\textbf{t} \in LTS \\ \rho(\textbf{t}) = p+1/2}}
    \hat{I}_{\textbf{t}}[F(\textbf{t})(X_{\cdot})]_{t_0,t} \, .
\end{alignat*}
%    \end{split}
%\end{equation}
}
%
%\\ \\

Taking into account the order of the symmetry group
$\sigma(\textbf{t})$ for a rooted tree $\textbf{t} \in TS$, we
finally obtain the relationship
\begin{equation}
    \begin{split}
    \sum_{\substack{\textbf{t} \in LTS \\ \rho(\textbf{t}) =
    1/2n}}
    \hat{I}_{\textbf{t}}[F(\textbf{t})(X_{\cdot})]_{t_0,t}
    = \sum_{\substack{\textbf{t} \in TS \\ \rho(\textbf{t}) =
    1/2 n}} \frac{I_{\textbf{t}}[F(\textbf{t})(X_{\cdot})]_{t_0,t}}{\sigma(\textbf{t})}
    \end{split}
\end{equation}
(see also \cite{Butcher03}) for all $n \in \mathbb{N}_0$ which
completes the proof. \hfill $\square$
%\\ \\

A similar stochastic Taylor expansion for the Stratonovich
SDE~(\ref{Intro-Ito-St-SDE1-integralform-Wm}) can be obtained
where the multiple stochastic integrals are defined with respect
to Stratonovich calculus.

\begin{Kor} \label{St-tree-expansion-exact-sol:Wm}
    For the solution process $(X_t)_{t \in [t_0,T]}$ of the
    Stratonovich
    SDE~(\ref{Intro-Ito-St-SDE1-integralform-Wm}) and for
    $p \in \tfrac{1}{2} \mathbb{N}_0$, $f :
    \mathbb{R}^d \rightarrow \mathbb{R}$
    with $f, a^i, b^{i,j}
    \in C^{2p+2}(\mathbb{R}^d, \mathbb{R})$
    for $i=1, \ldots, d$, $j=1, \ldots, m$, we obtain the
    expansion
    \begin{equation}
    \begin{split} \label{St-Tree-expansion-exact-sol-formula1:Wm}
        f(X_{t})
        = &\sum_{\substack{\textbf{t} \in TS \\ \rho(\textbf{t}) \leq
        p}} F(\textbf{t})(X_{t_0}) \frac{I_{\textbf{t};t_0,t}}{\sigma(\textbf{t})}
        + \underline{\mathcal{R}}_{p}(t,t_0)
    \end{split}
    \end{equation}
    $\Prob$-a.s.\ with remainder term
    \begin{equation}
    \label{St-Tree-expansion-exact-sol-remainder1:Wm}
        \begin{split}
        \underline{\mathcal{R}}_{p}(t,t_0) &=
        \sum_{\substack{\textbf{t} \in TS \\ \rho(\textbf{t}) = p+1/2}}
        \frac{I_{\textbf{t};t_0,t}[F(\textbf{t})(X_{\cdot})]}{\sigma(\textbf{t})}
%        \\
        + \sum_{\substack{\textbf{t} \in TS \\ \rho(\textbf{t}) = p}}
        \sum_{\textbf{u} \in H^0(\textbf{t})}
        \frac{I_{\textbf{t};t_0,t}[\int_{t_0}^{\cdot}
        %\sum_{k=1}^d a^k(X_s) \frac{\partial}{\partial x^k}
%        \hat{L}^0
        F(\textbf{u})(X_{s})
        \, {\mathrm{d}}s]}{\sigma(\textbf{t})}
        \end{split}
    \end{equation}
    provided all of the appearing multiple Stratonovich integrals exist.
\end{Kor}

We leave the proof of
Corollary~\ref{St-tree-expansion-exact-sol:Wm} to the reader since
it is analogously to that of
Theorem~\ref{Ito-tree-expansion-exact-sol:Wm} however with the
much simpler operator $\hat{L}^0$ instead of $L^0$.

Next, we give some results on the mean--square and mean
convergence of the obtained stochastic Taylor expansion by an
estimation for the remainder term. Therefore, we denote for $t \in
[t_0,T]$ and some $p \in \frac{1}{2} \mathbb{N}_0$ by
\begin{equation}
    Z_p(t) = \sum_{\substack{\textbf{t} \in TS \\ \rho(\textbf{t})
    \leq p}} F(\textbf{t})(X_{t_0})
    \frac{I_{\textbf{t};t_0,t}}{\sigma(\textbf{t})}
\end{equation}
the truncated stochastic Taylor expansion for
SDE~(\ref{Intro-Ito-St-SDE1-integralform-Wm}). Further, let $[p]$
denote the largest integer not exceeding $p$.

\begin{Sat} \label{Sat-Ito-error-estimeate}
    Let $X_{t_0} \in L^{2}(\Omega)$ and let $(X_t)_{t \in [t_0,T]}$ be the solution
    of the It{\^o} SDE~(\ref{Intro-Ito-St-SDE1-integralform-Wm}).
    Suppose that for
    all $\textbf{t} \in TS$ with $\rho(\textbf{t}) = p +
    \tfrac{1}{2}$ and for all $\textbf{t} \in H^0(\textbf{u}) \cup
    H^I(\textbf{u})$ with $\rho(\textbf{u})=p$ exists some constant
    $C>0$ such that
    \begin{equation} \label{Error-Vor1}
        \| F(\textbf{t})(x) \|^2 \leq C (1+ \|x\|^2) \, .
    \end{equation}
%    and for all $\textbf{t} \in TS$ with
%    $\rho(\textbf{t}) = p$
%    \begin{equation} \label{Error-Vor2}
%        \| L^0 F(\textbf{t})(x) \| \leq C (1+ \|x\|^2)^{1/2} \, .
%    \end{equation}
%    for some $C>0$.
%    For $t \in [t_0,T]$ let for the It{\^o}
%    SDE~(\ref{Intro-Ito-St-SDE1-integralform-Wm})
%    \begin{equation}
%        Z_p(t) = \sum_{\substack{\textbf{t} \in TS \\ \rho(\textbf{t}) \leq
%        p}} F(\textbf{t})(X_{t_0})
%        \frac{I_{\textbf{t};t_0,t}}{\sigma(\textbf{t})} \, .
%    \end{equation}
    Then there exists a constants $C_p>0$ depending on $p$ and
    a constant $C_L>0$
    depending on the Lipschitz constant of the drift and diffusion
    and on $T$, such that for all $t \in [t_0,T]$
    \begin{equation}
        \E (\| f(X_t) - Z_p(t) \|^{2}) \leq C_p (1+\E(\|X_{t_0}\|^2))
        \, \exp(C_L (t-t_0)) \,
        \frac{(t-t_0)^{2p+1}}{[p+\tfrac{1}{2}]!}
    \end{equation}
    for the mean--square truncation error and
    \begin{equation}
        \| \E (f(X_t) - Z_p(t)) \| \leq C_p
        (1+\E(\|X_{t_0}\|^2))^{1/2}
        \, \exp(C_L (t-t_0)) \,
        \frac{(t-t_0)^{p+\kappa}}{(p+\kappa)!}
    \end{equation}
    with $\kappa=1$ if $p \in \mathbb{N}_0$ and $\kappa = 1/2$ if $p \notin \mathbb{N}_0$
    for the mean truncation error.
\end{Sat}

{\bf Proof.}
%Apply Lemma~2.2 in \cite{Mil04} to the remainder
%$\mathcal{R}_{p}(t,t_0)$
Due to the Existence and Uniqueness Theorem \cite{Arn73,KS99}, for
$T>t_0$ there exists a constant $C>0$ which depends on
%$l \in \mathbb{N}$,
$T$ and the Lipschitz constants of $a$ and $b^j$, $j=1, \ldots,
m$, such that for all $t \in [t_0,T]$
\begin{equation}
    \E (\|X_t\|^{2}) \leq (1+ \E(\|X_{t_0}\|^{2}) \,
    \exp(C(t-t_0)) \, .
\end{equation}
Then, with (\ref{Error-Vor1})
%, (\ref{Error-Vor2})
and due to the It{\^o} isometry and with the Cauchy-Schwarz inequality
follows
\begin{equation*}
    \begin{split}
    &\E ( \| f(X_t) - Z_p(t)  \|^2 ) = \E ( \| \mathcal{R}_p(t,t_0) \|^2) \\
    %
%    &\leq \E (\| \sum_{\substack{\textbf{t} \in TS \\ \rho(\textbf{t}) = p+1/2}}
%    \frac{I_{\textbf{t};t_0,t}[F(\textbf{t})(X_{\cdot})]}{\sigma(\textbf{t})} \|^2)
%    + \E (\| \sum_{\substack{\textbf{t} \in TS \\ \rho(\textbf{t}) = p}}
%    \sum_{\textbf{u} \in H^0(\textbf{t})}
%    \frac{I_{\textbf{t};t_0,t}[\int_{t_0}^{\cdot} F(\textbf{u})(X_{s})
%    \, {\mathrm{d}}s]}{\sigma(\textbf{t})} \|^2) \\
%    &\quad + \E (\| \sum_{\substack{\textbf{t} \in TS \\ \rho(\textbf{t}) = p}}
%    \sum_{\textbf{u} \in H^I(\textbf{t})}
%    \frac{I_{\textbf{t};t_0,t}[\frac{1}{2} \int_{t_0}^{\cdot} F(\textbf{u})(X_{s})
%    \, {\mathrm{d}}s]}{\sigma(\textbf{t})} \|^2) \\
    %
    & \leq C_{1,p} \sum_{\substack{\textbf{t} \in TS \\ \rho(\textbf{t}) = p+1/2}}
    \E ( \| I_{\textbf{t};t_0,t}[F(\textbf{t})(X_{\cdot})] \|^2) \\
    &\quad + C_{2,p} \sum_{\substack{\textbf{t} \in TS \\ \rho(\textbf{t}) = p}}
    \sum_{\textbf{u} \in H^0(\textbf{t})}
    \E (\| I_{\textbf{t};t_0,t}[\int_{t_0}^{\cdot} F(\textbf{u})(X_{s})
    \, {\mathrm{d}}s] \|^2) \\
    &\quad + C_{3,p} \sum_{\substack{\textbf{t} \in TS \\ \rho(\textbf{t}) = p}}
    \sum_{\textbf{u} \in H^I(\textbf{t})}
    \E (\| I_{\textbf{t};t_0,t}[\frac{1}{2} \int_{t_0}^{\cdot} F(\textbf{u})(X_{s})
    \, {\mathrm{d}}s] \|^2) \\
    &\leq C_{4,p} \sum_{\substack{\textbf{t} \in TS \\ \rho(\textbf{t}) = p+1/2}}
    \frac{(t-t_0)^{2 \rho(\textbf{t})}}{[\rho(\textbf{t})]!}
    (1+\E(\|X_{t_0}\|^2)) \exp(C(t-t_0)) \\
    &\quad + C_{5,p} \sum_{\substack{\textbf{t} \in TS \\ \rho(\textbf{t}) = p}}
    \sum_{\textbf{u} \in H^0(\textbf{t}) \cup H^I(\textbf{t})}
    \frac{(t-t_0)^{2 \rho(\textbf{t})+2}}{[\rho(\textbf{t})+1]!}
    (1+\E(\|X_{t_0}\|^2)) \exp(C (t-t_0)) \\
    &\leq C_p (1+\E(\|X_{t_0}\|^2)) \exp(C (t-t_0)) \,
    \frac{(t-t_0)^{2p+1}}{[p+\tfrac{1}{2}]!} \, .
    \end{split}
\end{equation*}
Further, we obtain analogously with the Jensen inequality
\begin{equation*}
    \begin{split}
    &\| \E ( f(X_t) - Z_p(t) ) \|^2 =  \| \E(\mathcal{R}_p(t,t_0)) \|^2 \\
    & \leq C_{6,p} \sum_{\substack{\textbf{t} \in TS \\ \rho(\textbf{t}) = p+1/2}}
    \| \E ( I_{\textbf{t};t_0,t}[F(\textbf{t})(X_{\cdot})]) \|^2 \\
    &\quad + C_{7,p} \sum_{\substack{\textbf{t} \in TS \\ \rho(\textbf{t}) = p}}
    \sum_{\textbf{u} \in H^0(\textbf{t})}
    \| \E (I_{\textbf{t};t_0,t}[\int_{t_0}^{\cdot} F(\textbf{u})(X_{s})
    \, {\mathrm{d}}s]) \|^2 \\
    &\quad + C_{8,p} \sum_{\substack{\textbf{t} \in TS \\ \rho(\textbf{t}) = p}}
    \sum_{\textbf{u} \in H^I(\textbf{t})}
    \| \E (I_{\textbf{t};t_0,t}[\frac{1}{2} \int_{t_0}^{\cdot} F(\textbf{u})(X_{s})
    \, {\mathrm{d}}s]) \|^2 \\
    &\leq C_{9,p} \sum_{\substack{\textbf{t} \in TS \\ \rho(\textbf{t}) = p+1/2}}
    \frac{(t-t_0)^{2\rho(\textbf{t})+2 (\kappa-1/2)}}
    {((\rho(\textbf{t})+\kappa-\tfrac{1}{2})!)^2}
    (1+\E(\|X_{t_0}\|^2)) \exp(C (t-t_0)) \\
    &\quad + C_{10,p} \sum_{\substack{\textbf{t} \in TS \\ \rho(\textbf{t}) = p}}
    \sum_{\textbf{u} \in H^0(\textbf{t}) \cup H^I(\textbf{t})}
    \frac{(t-t_0)^{2 \rho(\textbf{t})+2}}{([\rho(\textbf{t})+1]!)^2}
    (1+\E(\|X_{t_0}\|^2)) \exp(C (t-t_0)) \\
    &\leq C_p (1+\E(\|X_{t_0}\|^2)) \exp(C (t-t_0)) \,
    \frac{(t-t_0)^{2p+2\kappa}}{((p+\kappa)!)^2} \, .
    \end{split}
\end{equation*}
\hfill $\square$

The remainder term $\underline{\mathcal{R}}_p(t,t_0)$ from the
Stratonovich expansion can be estimated analogously as in the
proof of Proposition~\ref{Sat-Ito-error-estimeate}. This follows
from Remark~5.2.8~\cite{KP99} because each multiple Stratonovich
integral can be written as a sum of multiple It{\^o} integrals of at
least the same mean and mean-square orders.

\begin{Kor}
    Let $X_{t_0} \in L^{2}(\Omega)$ and let $(X_t)_{t \in [t_0,T]}$
    be the solution of the Stra\-to\-no\-vich
    SDE~(\ref{Intro-Ito-St-SDE1-integralform-Wm}).
    Suppose that for
    all $\textbf{t} \in TS$ with $\rho(\textbf{t}) = p + \tfrac{1}{2}$ or
    $\rho(\textbf{t}) = p + 1$ exists some constant $C>0$ such
    that
%    and for all $\textbf{t} \in H^0(\textbf{u})$ with
%    $\rho(\textbf{u})=p$ and $\textbf{t} \in H^j(\textbf{v})$ with
%    $\rho(\textbf{v})=p+1/2$
    \begin{equation} \label{Error-Vor2}
        \| F(\textbf{t})(x) \|^2 \leq C (1+ \|x\|^2) \, .
    \end{equation}
%    and for all $\textbf{t} \in TS$ with
%    $\rho(\textbf{t}) = p$
%    \begin{equation} \label{Error-Vor2}
%        \| L^0 F(\textbf{t})(x) \| \leq C (1+ \|x\|^2)^{1/2} \, .
%    \end{equation}
%    for some $C>0$.
%    For $t \in [t_0,T]$ let for the Stratonovich
%    SDE~(\ref{Intro-Ito-St-SDE1-integralform-Wm})
%    \begin{equation}
%        Z_p(t) = \sum_{\substack{\textbf{t} \in TS \\ \rho(\textbf{t}) \leq
%        p}} F(\textbf{t})(X_{t_0})
%        \frac{I_{\textbf{t};t_0,t}}{\sigma(\textbf{t})} \, .
%    \end{equation}
    Then there exists a constant $C_p>0$ depending on $p$ and a constant $C_L>0$
    depending on the Lipschitz constant of the drift and diffusion
    and on $T$, such that for all $t \in [t_0,T]$
    \begin{equation}
        \E (\| f(X_t) - Z_p(t) \|^{2}) \leq C_p (1+\E(\|X_{t_0}\|^2))
        \, \exp(C_L (t-t_0)) \,
        \frac{(t-t_0)^{2p+1}}{[p+\tfrac{1}{2}]!}
    \end{equation}
    for the mean--square truncation error and
    \begin{equation}
        \| \E (f(X_t) - Z_p(t)) \| \leq C_p
        (1+\E(\|X_{t_0}\|^2))^{1/2}
        \, \exp(C_L (t-t_0)) \,
        \frac{(t-t_0)^{p+\kappa}}{(p+\kappa)!}
    \end{equation}
    with $\kappa=1$ if $p \in \mathbb{N}_0$ and $\kappa = 1/2$ if $p \notin \mathbb{N}_0$
    for the mean truncation error.
\end{Kor}
%
%
%
% ============================================================================
\section{Example}
\label{Sec:Example}
% ============================================================================
%
%
%
%
We consider SDE~(\ref{Intro-Ito-St-SDE1-integralform-Wm}) either
with respect to It{\^o} or Stratonovich calculus and give the
corresponding stochastic Taylor expansions
(\ref{Ito-Tree-expansion-exact-sol-formula1:Wm}) and
(\ref{St-Tree-expansion-exact-sol-formula1:Wm}) based on rooted
trees up to order $p=1.5$. Taking into account the trees presented
in Table~\ref{tabelle1} we obtain
{\allowdisplaybreaks
\begin{alignat*}{3}
%\begin{equation}
%    \begin{split}
    &f(X_t) = F(\textbf{t}_{0.1})(X_{t_0}) \,
    \frac{I_{\textbf{t}_{0.1};t_0,t}}{\sigma(\textbf{t}_{0.1})}
    + \sum_{1 \leq j_1 \leq m} F(\textbf{t}_{0.5.1})(X_{t_0}) \,
    \frac{I_{\textbf{t}_{0.5.1};t_0,t}}{\sigma(\textbf{t}_{0.5.1})}
    \\
    &+ F(\textbf{t}_{1.1})(X_{t_0}) \,
    \frac{I_{\textbf{t}_{1.1};t_0,t}}{\sigma(\textbf{t}_{1.1})}
    + \sum_{1 \leq j_1 \leq j_2 \leq m} F(\textbf{t}_{1.2})(X_{t_0}) \,
    \frac{I_{\textbf{t}_{1.2};t_0,t}}{\sigma(\textbf{t}_{1.2})} \\
    &+ \sum_{1 \leq j_1, j_2 \leq m} F(\textbf{t}_{1.3})(X_{t_0}) \,
    \frac{I_{\textbf{t}_{1.3};t_0,t}}{\sigma(\textbf{t}_{1.3})} +
    \sum_{1 \leq j_1 \leq m} F(\textbf{t}_{1.5.1})(X_{t_0}) \,
    \frac{I_{\textbf{t}_{1.5.1};t_0,t}}{\sigma(\textbf{t}_{1.5.1})} \\
    &+ \sum_{1 \leq j_1 \leq m} F(\textbf{t}_{1.5.2})(X_{t_0}) \,
    \frac{I_{\textbf{t}_{1.5.2};t_0,t}}{\sigma(\textbf{t}_{1.5.2})}
    + \sum_{1 \leq j_1 \leq m} F(\textbf{t}_{1.5.3})(X_{t_0}) \,
    \frac{I_{\textbf{t}_{1.5.3};t_0,t}}{\sigma(\textbf{t}_{1.5.3})}
    \\
    &+ \sum_{1 \leq j_1 \leq j_2 \leq j_3 \leq m} F(\textbf{t}_{1.5.4})(X_{t_0}) \,
    \frac{I_{\textbf{t}_{1.5.4};t_0,t}}{\sigma(\textbf{t}_{1.5.4})}
    + \sum_{1 \leq j_1, j_2, j_3 \leq m} F(\textbf{t}_{1.5.5})(X_{t_0}) \,
    \frac{I_{\textbf{t}_{1.5.5};t_0,t}}{\sigma(\textbf{t}_{1.5.5})}
    \\
    &+ \sum_{\substack{1 \leq j_1, j_2, j_3 \leq m \\ j_2 \leq j_3}}
    F(\textbf{t}_{1.5.6})(X_{t_0}) \,
    \frac{I_{\textbf{t}_{1.5.6};t_0,t}}{\sigma(\textbf{t}_{1.5.6})}
    + \sum_{1 \leq j_1, j_2, j_3 \leq m} F(\textbf{t}_{1.5.7})(X_{t_0}) \,
    \frac{I_{\textbf{t}_{1.5.7};t_0,t}}{\sigma(\textbf{t}_{1.5.7})} \\
    &+ R_{1.5}(t,t_0)
%    \end{split}
%\end{equation}
\end{alignat*}
}
with $R_{1.5}(t,t_0) = \mathcal{R}_{1.5}(t,t_0)$ in the case of
the It{\^o} SDE~(\ref{Intro-Ito-St-SDE1-integralform-Wm}) and
$R_{1.5}(t,t_0) = \underline{\mathcal{R}}_{1.5}(t,t_0)$ in the
case of the Stratonovich
SDE~(\ref{Intro-Ito-St-SDE1-integralform-Wm}). Now, we apply the
indicator function $1_{\{* \neq \circ\}}$ which vanishes in the
case of Stratonovich calculus. Calculating the elementary
differentials and the multiple integrals, we get
{\allowdisplaybreaks
%\begin{equation}
%    \begin{split}
\begin{alignat*}{3}
    f(X_t) &= f(X_{t_0})
    + \sum_{1 \leq j_1 \leq m} \sum_{J=1}^d (\frac{\partial f}{\partial
    x^J} \, b^{J,j_1}) (X_{t_0}) \,
    \mathcal{I}_{(j_1)}[1]_{t_0,t} \\
    %\frac{I_{\textbf{t}_{0.5.1};t_0,t}}{\sigma(\textbf{t}_{0.5.1})} \\
    &+ \sum_{J=1}^d (\frac{\partial f}{\partial
    x^J} \, a^{J})(X_{t_0}) \,
    \mathcal{I}_{(0)}[1]_{t_0,t} \\
    %\frac{I_{\textbf{t}_{1.1};t_0,t}}{\sigma(\textbf{t}_{1.1})} \\
    &+ \sum_{1 \leq j_1 \leq j_2 \leq m} \sum_{J,K=1}^d (\frac{\partial^2 f}
    {\partial x^J \partial x^K} \, b^{J,j_1} \, b^{K,j_2})(X_{t_0}) \\
    & \times \frac{\mathcal{I}_{(j_1,j_2)}[1]_{t_0,t} + \mathcal{I}_{(j_1,j_2)}[1]_{t_0,t}
    + \mathcal{I}_{(0)}[1_{\{j_1=j_2 \wedge * \neq
    \circ\}}]_{t_0,t}}{1_{\{j_1 \neq j_2\}}+2 \cdot 1_{\{j_1=j_2\}}} \\
    %\frac{I_{\textbf{t}_{1.2};t_0,t}}{\sigma(\textbf{t}_{1.2})} \\
    &+ \sum_{1 \leq j_1, j_2 \leq m} \sum_{J,K=1}^d (\frac{\partial f}{\partial
    x^J} \, \frac{\partial b^{J,j_1}}{x^K} \, b^{K,j_2})(X_{t_0}) \,
    \mathcal{I}_{(j_2,j_1)}[1]_{t_0,t} \\
    %\frac{I_{\textbf{t}_{1.3};t_0,t}}{\sigma(\textbf{t}_{1.3})} \\
    &+ \sum_{1 \leq j_1 \leq m} \sum_{J,K=1}^d (\frac{\partial f}{\partial
    x^J} \, \frac{\partial a^J}{\partial x^K} \, b^{K,j_1})(X_{t_0}) \,
    \mathcal{I}_{(j_1,0)}[1]_{t_0,t} \\
    %\frac{I_{\textbf{t}_{1.5.1};t_0,t}}{\sigma(\textbf{t}_{1.5.1})} \\
    &+ \sum_{1 \leq j_1 \leq m} \sum_{J,K=1}^d (\frac{\partial f}{\partial
    x^J} \, \frac{\partial b^{J,j_1}}{\partial x^K} \, a^K)(X_{t_0}) \,
    \mathcal{I}_{(0,j_1)}[1]_{t_0,t} \\
    %\frac{I_{\textbf{t}_{1.5.2};t_0,t}}{\sigma(\textbf{t}_{1.5.2})} \\
    &+ \sum_{1 \leq j_1 \leq m} \sum_{J,K=1}^d (\frac{\partial^2 f}{\partial
    x^J \partial x^K} \, a^J \, b^{K,j_1})(X_{t_0}) \,
    (\mathcal{I}_{(j_1,0)}[1]_{t_0,t} +
    \mathcal{I}_{(0,j_1)}[1]_{t_0,t}) \\
    %\frac{I_{\textbf{t}_{1.5.3};t_0,t}}{\sigma(\textbf{t}_{1.5.3})} \\
    &+ \sum_{1 \leq j_1 \leq j_2 \leq j_3 \leq m} \sum_{J,K,L=1}^d (\frac{\partial^3 f}{\partial
    x^J \partial x^K \partial x^L} \, b^{J,j_1} \, b^{K,j_2} \, b^{L,j_3})(X_{t_0}) \\ %\,
    &\times (\mathcal{I}_{(j_3,j_2,j_1)}[1]_{t_0,t} + \mathcal{I}_{(j_3,j_1,j_2)}[1]_{t_0,t}
    + \mathcal{I}_{(j_1,j_3,j_2)}[1]_{t_0,t} +
    \mathcal{I}_{(j_2,j_3,j_1)}[1]_{t_0,t} \\
    &+ \mathcal{I}_{(j_2,j_1,j_3)}[1]_{t_0,t} + \mathcal{I}_{(j_1,j_2,j_3)}[1]_{t_0,t}
    + \mathcal{I}_{(0,j_2)}[1_{\{j_1=j_3 \wedge * \neq \circ\}}]_{t_0,t} \\
    &+ \mathcal{I}_{(j_3,0)}[1_{\{j_1=j_2 \wedge * \neq \circ\}}]_{t_0,t}
    + \mathcal{I}_{(0,j_3)}[1_{\{j_1=j_2 \wedge * \neq \circ\}}]_{t_0,t}
    + \mathcal{I}_{(j_2,0)}[1_{\{j_1=j_3 \wedge * \neq \circ\}}]_{t_0,t} \\
    &+ \mathcal{I}_{(0,j_1)}[1_{\{j_2=j_3 \wedge * \neq \circ\}}]_{t_0,t}
    + \mathcal{I}_{(j_1,0)}[1_{\{j_2=j_3 \wedge * \neq \circ\}}]_{t_0,t})
    \\
    &\times \frac{1}{1_{\{j_1\neq j_2 \neq j_3 \neq j_1 \}}
    + 2 \cdot 1_{\{j_1=j_2\neq j_3 \vee j_1=j_3 \neq j_2 \vee j_2=j_3
    \neq j_1\}} + 6 \cdot 1_{\{j_1=j_2=j_3\}} } \\
    %\frac{I_{\textbf{t}_{1.5.4};t_0,t}}{\sigma(\textbf{t}_{1.5.4})} \\
    &+ \sum_{1 \leq j_1, j_2, j_3 \leq m} \sum_{J,K,L=1}^d (\frac{\partial^2 f}{\partial
    x^J \partial x^K} \, \frac{\partial b^{J,j_1}}{\partial x^L} \, b^{L,j_2} \, b^{K,j_3})
    (X_{t_0}) \\ %\,
    &\times (\mathcal{I}_{(j_3,j_2,j_1)}[1]_{t_0,t} + \mathcal{I}_{(j_2,j_3,j_1)}[1]_{t_0,t}
    + \mathcal{I}_{(j_2,j_1,j_3)}[1]_{t_0,t} \\
    &+ \mathcal{I}_{(0,j_1)}[1_{\{j_2=j_3 \wedge * \neq \circ\}}]_{t_0,t}
    + \mathcal{I}_{(j_2,0)}[1_{\{j_1=j_3 \wedge * \neq \circ\}}]_{t_0,t}) \\
    %\frac{I_{\textbf{t}_{1.5.5};t_0,t}}{\sigma(\textbf{t}_{1.5.5})} \\
    &+ \sum_{\substack{1 \leq j_1, j_2, j_3 \leq m \\ j_2 \leq j_3}}
    \sum_{J,K,L=1}^d (\frac{\partial f}{\partial
    x^J} \, \frac{\partial^2 b^{J,j_1}}{\partial x^K \partial x^L} \, b^{K,j_2} \, b^{L,j_3})
    (X_{t_0}) \\ %\,
    &\times \frac{\mathcal{I}_{(j_3,j_2,j_1)}[1]_{t_0,t} + \mathcal{I}_{(j_2,j_3,j_1)}[1]_{t_0,t}
    + \mathcal{I}_{(0,j_1)}[1_{\{j_2=j_3 \wedge * \neq
    \circ\}}]_{t_0,t}}{1_{\{j_2 \neq j_3\}} + 2 \cdot 1_{\{j_2=j_3\}} } \\
    %\frac{I_{\textbf{t}_{1.5.6};t_0,t}}{\sigma(\textbf{t}_{1.5.6})} \\
    &+ \sum_{1 \leq j_1, j_2, j_3 \leq m} \sum_{J,K,L=1}^d (\frac{\partial f}{\partial
    x^J} \, \frac{\partial b^{J,j_1}}{\partial x^K} \, \frac{\partial b^{K,j_2}}{\partial x^L}
    \, b^{L,j_3})(X_{t_0}) \,
    \mathcal{I}_{(j_3,j_2,j_1)}[1]_{t_0,t} \\
    %\frac{I_{\textbf{t}_{1.5.7};t_0,t}}{\sigma(\textbf{t}_{1.5.7})}
    &+ R_{1.5}(t,t_0) \, .
\end{alignat*}
%    \end{split}
%\end{equation}
}
As the main advantage, each elementary differential corresponds
exactly to one tree. Further, all equal elementary differentials
are pooled with a corresponding weight which may be a sum of some
multiple stochastic integrals. For example, if we apply the
calculated expansion of order $p=1.5$ to the It{\^o} SDE
\begin{equation}
    {\mathrm{d}} X_t = \mu \, X_t \, {\mathrm{d}}t + \sigma \, X_t \,
    {\mathrm{d}}W_t, \qquad X_{t_0}=x_0,
\end{equation}
with $d=m=1$ and $\mu,\sigma \in \mathbb{R}$, then we obtain for
$f(x)=x$ the expansion
\begin{equation}
    \begin{split}
    X_t &= x_0 + \sigma \, x_0 \, \mathcal{I}_{(1)} + \mu \, x_0 \,
    \mathcal{I}_{(0)} + \sigma^2 \, x_0 \, \mathcal{I}_{(1,1)} + \mu \,
    \sigma \, x_0 \, \mathcal{I}_{(1,0)} \\
    &+ \sigma \, \mu \, x_0 \,
    \mathcal{I}_{(0,1)} + \sigma^3 \, x_0 \, \mathcal{I}_{(1,1,1)} +
    \mathcal{R}_{1.5}(t,t_0) \, .
    \end{split}
\end{equation}
Here, we obtain with Proposition~\ref{Sat-Ito-error-estimeate}
that $(E(\| \mathcal{R}_{1.5}(t,t_0) \|^2))^{1/2} =
O((t-t_0)^{2})$ as an estimate for the mean--square error and $\|
E (\mathcal{R}_{1.5}(t,t_0)) \| = O((t-t_0)^{2})$ for the mean
error.
\begin{longtable}[tbp]{|c|c||c|c||c|c|}
\caption{Trees $\textbf{t} \in TS$ of order $\rho(\textbf{t}) \leq
2$ with $j_1, j_2, j_3, j_4 \in \{1,\ldots,m\}$.} \label{tabelle1}
\\
    \hline
    $\textbf{t}$ & tree & $\textbf{t}$ &
    tree & $\textbf{t}$ & tree \\
    \hline
    \hline
    \endfirsthead
    \hline
    $\textbf{t}$ & tree & $\textbf{t}$ &
    tree & $\textbf{t}$ & tree \\
    \hline
    \hline
    \endhead
%
%    \hline
    \hline
    \endfoot
%
%    \hline
    \hline
    \endlastfoot
    $\textbf{t}_{0.1}$ & $\gamma$ & $\textbf{t}_{0.5.1}$ & $[\tau_{j_1}]_{\gamma}$ & & \\
%    \cline{1-8}
    \hline
    $\textbf{t}_{1.1}$ & $[\tau_0]_{\gamma}$ & $\textbf{t}_{1.2}$ & $[\tau_{j_1},\tau_{j_2}]_{\gamma}$
    & $\textbf{t}_{1.3}$ & $[[\tau_{j_2}]_{j_1}]_{\gamma}$ \\
%    \cline{1-8}
    \hline
    $\textbf{t}_{1.5.1}$ & $[[\tau_{j_1}]_0]_{\gamma}$ & $\textbf{t}_{1.5.2}$ & $[[\tau_0]_{j_1}]_{\gamma}$
    & $\textbf{t}_{1.5.3}$ & $[\tau_0,\tau_{j_1}]_{\gamma}$ \\
    $\textbf{t}_{1.5.4}$ & $[\tau_{j_1},\tau_{j_2},\tau_{j_3}]_{\gamma}$
    & $\textbf{t}_{1.5.5}$ & $[[\tau_{j_2}]_{j_1},\tau_{j_3}]_{\gamma}$
    & $\textbf{t}_{1.5.6}$ & $[[\tau_{j_2},\tau_{j_3}]_{j_1}]_{\gamma}$ \\
    $\textbf{t}_{1.5.7}$ & $[[[\tau_{j_3}]_{j_2}]_{j_1}]_{\gamma}$ & & & & \\
%    \cline{1-8}
    \hline
    $\textbf{t}_{2.1}$ & $[[\tau_0]_0]_{\gamma}$ & $\textbf{t}_{2.2}$ & $[\tau_0,\tau_0]_{\gamma}$
    & $\textbf{t}_{2.3}$ & $[[[\tau_{j_2}]_{j_1}]_0]_{\gamma}$ \\
    $\textbf{t}_{2.4}$ & $[[\tau_{j_1},\tau_{j_2}]_0]_{\gamma}$
    & $\textbf{t}_{2.5}$ & $[\tau_{j_1},[\tau_{j_2}]_0]_{\gamma}$
    & $\textbf{t}_{2.6}$ & $[[\tau_{j_2}]_{j_1},\tau_0]_{\gamma}$ \\
    $\textbf{t}_{2.7}$ & $[\tau_{j_1},\tau_{j_2},\tau_0]_{\gamma}$
    & $\textbf{t}_{2.8}$ & $[\tau_{j_1},[\tau_0]_{j_2}]_{\gamma}$
    & $\textbf{t}_{2.9}$ & $[[[\tau_0]_{j_2}]_{j_1}]_{\gamma}$ \\
    $\textbf{t}_{2.10}$ & $[[\tau_{j_2},\tau_0]_{j_1}]_{\gamma}$
    & $\textbf{t}_{2.11}$ & $[\tau_{j_1},\tau_{j_2},\tau_{j_3},\tau_{j_4}]_{\gamma}$
    & $\textbf{t}_{2.12}$ & $[\tau_{j_1},\tau_{j_2},[\tau_{j_4}]_{j_3}]_{\gamma}$ \\
    $\textbf{t}_{2.13}$ & $[\tau_{j_1},[\tau_{j_3},\tau_{j_4}]_{j_2}]_{\gamma}$
    & $\textbf{t}_{2.14}$ & $[\tau_{j_1},[[\tau_{j_4}]_{j_3}]_{j_2}]_{\gamma}$
    & $\textbf{t}_{2.15}$ & $[[\tau_{j_2}]_{j_1},[\tau_{j_4}]_{j_3}]_{\gamma}$ \\
    $\textbf{t}_{2.16}$ & $[[\tau_{j_2},\tau_{j_3},\tau_{j_4}]_{j_1}]_{\gamma}$
    & $\textbf{t}_{2.17}$ & $[[\tau_{j_2},[\tau_{j_4}]_{j_3}]_{j_1}]_{\gamma}$
    & $\textbf{t}_{2.18}$ & $[[[\tau_{j_3},\tau_{j_4}]_{j_2}]_{j_1}]_{\gamma}$ \\
    $\textbf{t}_{2.19}$ & $[[[[\tau_{j_4}]_{j_3}]_{j_2}]_{j_1}]_{\gamma}$
    & $\textbf{t}_{2.20}$ & $[[[\tau_{j_2}]_0]_{j_1}]_{\gamma}$ & & \\
%    \cline{1-8}
    \hline
\end{longtable}

\end{document}